\documentclass[11pt]{article}

\usepackage{amsmath,amssymb,enumitem}
\usepackage{pstricks,pst-node,xy}
\xyoption{all}

\usepackage{xr}
\numberwithin{equation}{section}

\newtheorem{thm}{Theorem}
\newtheorem{thmv}{Theorem}

\newtheorem{prp}{Proposition}[section]
\newtheorem{lmm}[prp]{Lemma}

\def\lra{\longrightarrow}
\def\Lra{\Longrightarrow}

\def\lan{\langle}
\def\ran{\rangle}

\def\e_ref#1{(\ref{#1})}
\def\ov#1{\overline{#1}}

\def\ti#1{\tilde{#1}}
\def\wt#1{\widetilde{#1}}
\def\smsize#1{\begin{small}#1\end{small}}
\def\blr#1{\big\lan{#1}\big\ran}
\def\lr#1{\lan{#1}\ran}
\def\Blr#1{\Big\lan{#1}\Big\ran}

\def\al{\alpha}
\def\be{\beta}
\def\ga{\gamma}
\def\io{\iota}

\def\la{\lambda}
\def\om{\omega}
\def\si{\sigma}

\def\vph{\varphi}
\def\vr{\varrho}

\def\De{\Delta}
\def\Ga{\Gamma}
\def\La{\Lambda}
\def\Si{\Sigma}
\def\Th{\Theta}

\def\ale{\aleph}
\def\bpar{\bar\partial}
\def\eset{\emptyset}
\def\i{\infty}

\def\A{\mathcal A}
\def\B{\mathcal B}
\def\c{\mathbf c}
\def\C{\mathbb C}
\def\cC{\mathcal C}

\def\cD{\mathcal D}
\def\E{\mathbb E}
\def\G{\mathfrak G}

\def\J{\mathcal J}

\def\M{\mathfrak M}
\def\cM{\mathcal M}

\def\O{\mathcal O}
\def\P{\mathbb P}
\def\Q{\mathbb Q}
\def\R{\mathbb R}
\def\T{\mathcal T}
\def\U{\mathfrak U}
\def\V{\mathcal V}
\def\X{\mathfrak X}
\def\Z{\mathbb Z}
\def\cZ{\mathcal Z}

\def\deg{\textnormal{deg\,}}
\def\dim{\textnormal{dim}}
\def\es{\textnormal{es}}
\def\ev{\textnormal{ev}}
\def\gd{\textnormal{gd}}
\def\GW{\textnormal{GW}}
\def\Hom{\textnormal{Hom}}
\def\id{\textnormal{id}}
\def\rk{\textnormal{rk}}

\begin{document}

\title{Standard vs.~Reduced\\ 
Genus-One Gromov-Witten Invariants}
\author{Aleksey Zinger\thanks{Partially supported by a Sloan fellowship 
and DMS Grant 0604874}}
\date{\today}

\maketitle

\begin{abstract}
\noindent
We give an explicit formula for the difference between the standard and reduced
genus-one Gromov-Witten invariants.
Combined with previous work on geometric properties of the latter,
this paper makes it possible to compute the standard genus-one GW-invariants
of complete intersections.
In particular, we obtain a closed formula for the genus-one GW-invariants
of a Calabi-Yau projective hypersurface and verify a recent mirror symmetry
prediction for a sextic fourfold as a special case.
\end{abstract}

\tableofcontents

\section{Introduction}
\label{intro_sec}

\noindent
Gromov-Witten invariants are counts of holomorphic curves in symplectic manifolds
that play prominent roles in string theory, symplectic topology, and algebraic geometry.
A variety of predictions concerning these invariants have arisen from string theory,
only some of which have been verified mathematically.
These invariants are generally difficult to compute, especially in positive genera.
For example, the 1991 prediction of~\cite{CaDGP} for the genus-zero GW-invariants 
of the quintic threefold was mathematically confirmed in the mid-1990s,
while even low-degree cases of the 1993 prediction of \cite{BCOV} for its 
genus-one GW-invariants remained unaccessible for another seven years.\\

\noindent
In contrast to the genus-zero case, the expected hyperplane (or hyperplane-section)
relation between  GW-invariants of a complete intersection and 
those of the ambient space do not hold in positive genera;
see Subsection~\ref{bcov1-appr_subs} in~\cite{bcov1} for more details.
This issue is entirely avoided in~\cite{Ga} and~\cite{MP1} by approaching  
GW-invariants of complete intersections through degeneration techniques.
The methods of~\cite{Ga} and~\cite{MP1} can be used for low-degree
checks of~\cite{BCOV},
but they do not seem to provide a ready platform for an application 
of combinatorial techniques as in the genus-zero case.
In contrast, the failure of the expected hyperplane relation for genus-one invariants
is addressed in~\cite{g1comp2} and~\cite{LZ} by defining {\tt reduced} 
genus-one GW-invariants and showing that these invariants do satisfy 
the hyperplane relation, respectively.
Combined with~\cite{VaZ}, this last approach provides a platform suitable for 
an application of combinatorial techniques (the classical localization theorem) 
and has led to a closed formula for the reduced GW-invariants of 
a Calabi-Yau projective hypersurface; see~\cite{bcov1}.
For symplectic manifolds of real dimension~six, the standard and 
reduced genus-one GW-invariants without descendants differ by a multiple of 
the genus-zero  GW-invariant.\footnote{For the purposes of this statement, 
all three  GW-invariants are viewed as linear functionals on the same vector space 
(consisting of tuples of homology elements in the manifold).}
The prediction of~\cite{BCOV} is thus fully verified in~\cite{bcov1} 
as a special case.\\

\noindent
{\it Remark:} The ranges of applicability of the two degeneration methods 
and of the reduced-invariants method are very different, with a 
relatively limited overlap.
The wider-ranging degeneration method, that of~\cite{MP1}, 
can in principle be used to compute arbitrary-genus GW-invariants 
of low-degree low-dimension complete intersections.
While the computation in each case is generally difficult, 
this method has been used to compute genus-one and genus-two GW-invariants of
the Enriques Calabi-Yau threefold; see~\cite{MP2}.
On the other hand, the reduced-invariants method applies to arbitrary complete intersections,
but at this point in the genus-one case only.\\

\noindent
In theory, reduced genus-one GW-invariants are not new invariants,
as by Proposition~\ref{g1comp2-bdcontr_prp} in~\cite{g1comp2}
the difference between these invariants and the standard ones
is some combination of genus-zero GW-invariants.
However, this is not part of the definition of reduced genus-one invariants 
and the relation described by Proposition~\ref{g1comp2-bdcontr_prp} in~\cite{g1comp2}
is not convenient for immediate applications; see Proposition~\ref{bdcontr_prp} below.
In this paper, we determine an explicit relation; see 
Theorems~\ref{main_thm1} and~\ref{main_thm2}.
Combining this relation with the closed formula for the reduced genus-one GW-invariants
of a Calabi-Yau hypersurface derived in~\cite{bcov1}, 
we then obtain a closed formula for the standard genus-one GW-invariants
of a Calabi-Yau hypersurface; see Theorem~\ref{cy_thm}.
The mirror symmetry prediction of~\cite{KP} for a sextic fourfold is confirmed 
by the $n\!=\!6$ case of Theorem~\ref{cy_thm}.
The $n\!>\!6$ cases of Theorem~\ref{cy_thm} go beyond even predictions,
as far as the author is aware.\\

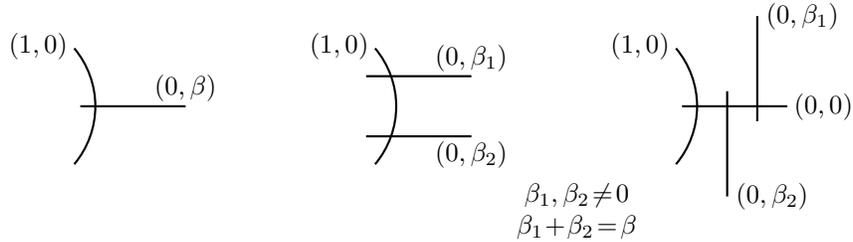
\begin{figure}
\begin{pspicture}(-1.4,-1.8)(10,2)
\psset{unit=.4cm}
\psarc(3,0){3}{-40}{40}\rput(4.1,2){\smsize{$(1,0)$}}
\psline(5.5,0)(9,0)\rput(9,.6){\smsize{$(0,\be)$}}
\psarc(13,0){3}{-40}{40}\rput(14.1,2){\smsize{$(1,0)$}}
\psline(15,1)(18.5,1)\rput(18.5,1.6){\smsize{$(0,\be_1)$}}
\psline(15,-1)(18.5,-1)\rput(18.5,-1.6){\smsize{$(0,\be_2)$}}
\rput(22,-3.5){\smsize{\begin{tabular}{c}$\be_1,\be_2\!\neq\!0$\\
$\be_1\!+\!\be_2\!=\!\be$\end{tabular}}}
\psarc(23,0){3}{-40}{40}\rput(24.1,2){\smsize{$(1,0)$}}
\psline(25.5,0)(29,0)\rput(30.2,0){\smsize{$(0,0)$}}
\psline(27,.5)(27,-3)\rput(28.5,-3){\smsize{$(0,\be_2)$}}
\psline(28,-.5)(28,3)\rput(29.5,3){\smsize{$(0,\be_1)$}}
\end{pspicture}
\caption{Generic elements of $\ov\M_{1,k}^{1,\eset}(X,\be)$, $\ov\M_{1,k}^{2,\eset}(X,\be)$, 
and $\ov\M_{1,k}^{1,\eset}(X,\be)\!\cap\!\ov\M_{1,k}^{2,\eset}(X,\be)$, respectively;
the lines and curves represent the components of the domain of a stable map and 
the pair of integers next to each component indicates the genus of the component
and the degree of the map on the component.}
\label{intro_fig}
\end{figure}

\noindent
It is interesting to observe that only one boundary stratum of a partially regularized
moduli space of genus-one stable maps accounts for the difference between
the standard and reduced genus-one GW-invariants without descendants.
This implies that if $X$ is a sufficiently regular almost Kahler manifold 
(e.g.~a low-degree projective hypersurface), only two strata of the moduli space
$\ov\M_{1,k}(X,\be)$ of degree-$\be$ genus-one stable maps to $X$ with $k$ marked points
contribute to the genus-one GW-invariants without descendants:
\begin{enumerate}[label=(\roman*)]
\item the main stratum $\M_{1,k}^0(X,\be)$ consisting of stable maps 
from smooth domains;
\item the stratum $\M_{1,k}^{1,\eset}(X,\be)$ consisting of stable maps 
from a union of a smooth genus-one curve and a~$\P^1$, with the map constant 
on the genus-one curve and with all $k$ marked points lying on the latter; 
see the first diagram in Figure~\ref{intro_fig}.
\end{enumerate}
At a first glance, this statement may seem to contradict reality,
as well as Theorems~\ref{main_thm1} and~\ref{main_thm2}.
For example, if $n\!\ge\!4$, the formula for the difference in 
Theorems~\ref{main_thm1} and~\ref{main_thm2} involves the GW-invariant $\GW_{(2,\eset)}^\be$
that counts two-component rational curves; this is consistent with~\cite{KP}.
This term may appear to arise from the stratum $\M_{1,k}^{2,\eset}(X,\be)$
of $\ov\M_{1,k}(X,\be)$ consisting of maps from a smooth genus-one curve~$\cC_P$
with two $\P^1$'s attached directly to~$\cC_P$ so that the map is constant on~$\cC_P$;
see the middle diagram in Figure~\ref{intro_fig}.
In fact, $\GW_{(2,\eset)}^\be$ arises from a homology class on 
$\ov\M_{0,k+1}(X,\be)$, or equivalently from the closure of a boundary stratum
of $\ov\M_{1,k}^{1,\eset}(X,\be)$.
This boundary stratum is the intersection of $\ov\M_{1,k}^{1,\eset}(X,\be)$
with $\M_{1,k}^{2,\eset}(X,\be)$; see the last diagram in Figure~\ref{intro_fig}.\\

\noindent
After setting up notation for GW-invariants in Subsection~\ref{gwinv_subs},
we state the main theorem of this paper is Subsection~\ref{mainthm_subs}.
Theorem~\ref{main_thm1} expresses the difference between the standard and
reduced genus-one GW-invariants as a linear combination of genus-zero invariants.
The coefficients in this linear combination are top intersections of tautological
classes on the blowups of moduli spaces of genus-one curves constructed
in Subsection~\ref{g1desing-curve1bl_subs} in~\cite{VaZ}
and reviewed in Subsection~\ref{curvbl_subs} below.
These are computable through the recursions obtained in~\cite{g1desing2}
and restated in Subsection~\ref{mainthm_subs} below;
\e_ref{psiclass_e} gives an explicit formula for these coefficients  
when no descendants are involved.
We then deduce a more compact version of Theorem~\ref{main_thm1};
Theorem~\ref{main_thm2} involves certain (un-)twisted $\psi$-classes
and coefficients that satisfy simpler recursions than the coefficients
in Theorem~\ref{main_thm1}.
The descendant-free case of Theorem~\ref{main_thm2}, \e_ref{main_red_e},
is used in Subsection~\ref{cy_subs} to obtain a closed formula 
for the difference between the two genus-one GW-invariants of
a Calabi-Yau hypersurface from a closed formula for its genus-zero
GW-invariants obtained in~\cite{MirSym}; see Lemma~\ref{cydiff_lmm}.
Using a closed formula for the reduced genus-one GW-invariants of
such a hypersurface derived in~\cite{bcov1}, we thus obtain a closed 
formula for its standard genus-one GW-invariants.\\

\noindent
Theorem~\ref{main_thm1} is proved in Section~\ref{pf_sec}.
In addition to reviewing the blowup construction of
Subsection~\ref{g1desing-curve1bl_subs} in~\cite{VaZ},
Subsection~\ref{curvbl_subs} describes natural bundle homomorphisms
over moduli spaces of genus-one curves and their twisted versions.
These are used to describe the difference between the two genus-one 
GW-invariants in Subsection~\ref{diff_subs} and to compute it
explicitly in Subsection~\ref{comp_subs}, respectively.
The blowup construction of Subsection~\ref{g1desing-curve0bl_subs} in~\cite{VaZ}
for moduli spaces of genus-{\it{zero}} curves is used in 
Subsection~\ref{curvbl0_subs} to obtain a formula for top intersections of
tautological classes on blowups of moduli spaces of genus-{\it{one}} curves;
Proposition~\ref{psiform_prp} is used at the end of Subsection~\ref{comp_subs}.
Subsection~\ref{diff_subs} reviews the relevant aspects of~\cite{g1comp2},
concluding with a description of the difference between the two genus-one
GW-invariants; see Proposition~\ref{bdcontr_prp}.
This difference can be computed explicitly through the direct,
but rather involved, setup of~\cite{g1}.
Subsection~\ref{comp_subs} instead presents a more conceptual approach
motivated by the blowup construction of 
Section~\ref{g1desing-map0bl_sec} in~\cite{VaZ} 
for moduli spaces of genus-zero maps.\\

\noindent
The author would like to thank J.~Li for first drawing the author's attention to
the issue of computing genus-one GW-invariants of projective hypersurfaces and
R.~Pandharipande for pointing out the mirror symmetry prediction for a 
sextic fourfold in~\cite{KP}.

\section{Overview}
\label{overview_sec}

\subsection{Gromov-Witten Invariants}
\label{gwinv_subs}

\noindent
Let $(X,\om,\J)$ be a compact symplectic manifold with a compatible almost complex structure.
If $g,k\!\in\!\bar\Z^+$ are nonnegative integers and $\be\!\in\!H_2(X;\Z)$, we denote by 
$\ov\M_{g,k}(X,\be;\J)$ the moduli space of genus-$g$ degree-$\be$ $\J$-holomorphic maps
into $X$ with $k$ marked points.
For each $j\!=\!1,\ldots,k$, let
$$\ev_j\!:\ov\M_{g,k}(X,\be;\J)\lra X$$
be the evaluation map at the $j$th marked point and let 
$$\psi_j\in H^2\big(\ov\M_{g,k}(X,\be;\J)\big)$$
be the first chern class of the universal cotangent line bundle at the $j$th marked point.
More generally, if $J$ is a finite set, we denote by $\ov\M_{0,J}(X,\be;\J)$ the moduli space of genus-$0$ degree-$\be$ $\J$-holomorphic maps into $X$ with marked points indexed by the set~$J$
and by
$$\ev_j\!:\ov\M_{0,J}(X,\be;\J)\lra X, \quad
\psi_j\in H^2\big(\ov\M_{0,J}(X,\be;\J)\big), \qquad\forall\,j\in J,$$
the corresponding evaluation maps and $\psi$-classes.
If $\be\!\neq\!0$, for each $J'\!\subset\!J$ there is a well-defined forgetful map
$$f_{J'}\!:\ov\M_{0,J}(X,\be;\J)\lra\ov\M_{0,J-J'}(X,\be;\J),$$
obtained by dropping the marked points indexed by $J'$ from the domain of
every stable map in $\ov\M_{0,J}(X,\be;\J)$ and contracting the unstable components
of the resulting map.
Let
$$\ti\psi_j\equiv f_{J-j}^*\psi_j\in H^2\big(\ov\M_{0,J}(X,\be;\J)\big)$$
be {\tt the untwisted $j$th $\psi$-class}.\\

\noindent
We also define moduli spaces of tuples of genus-zero stable maps.
If $m\!\in\!\bar\Z^+$, let
$$[m]=\big\{i\!\in\!\Z^+\!:1\!\le\!i\!\le\!m\big\}.$$
If $m\!\in\!\Z^+$ and $J$ is a finite set, we define
\begin{equation*}\begin{split}
\ov\M_{(m,J)}(X,\be;\J)&= 
\bigg\{(b_i)_{i\in[m]}\in\prod_{i=1}^{i=m}\ov\M_{0,\{0\}\sqcup J_i}(X,\be_i;\J)\!:
\be_i\!\in\!H_2(X;\Z)\!-\!\{0\},~J_i\!\subset\!J;\\
&\hspace{1.2in}~~
\sum_{i=1}^{i=m}\be_i\!=\!\be,~\bigsqcup_{i=1}^{i=m}J_i\!=\!J,~
\ev_0(b_i)\!=\!\ev_0(b_{i'})~\forall\, i,i'\!\in\![m]\bigg\}.
\end{split}\end{equation*}
There is a well-defined evaluation map
$$\ev_0\!: \ov\M_{(m,J)}(X,\be;\J)\lra X, \qquad
(b_i)_{i\in[m]}\lra \ev_0(b_i),$$
where $i$ is any element of $[m]$.
For each $i\!\in\![m]$, let 
$$\pi_i\!: \ov\M_{(m,J)}(X,\be;\J)\lra 
\bigsqcup_{\be_i\in H_2(X;\Z)}\bigsqcup_{J_i\subset J}
\ov\M_{0,\{0\}\sqcup J_i}(X,\be_i;\J)$$
be the projection onto the $i$th component.
If $p\!\in\!\bar\Z^+$, we define
$$\eta_p,\ti\eta_p\in H^{2p}\big(\ov\M_{(m,J)}(X,\be;\J)\big)$$
to be the sums of all degree-$p$ monomials in 
$$\big\{\pi_i^*\psi_0\!:~i\!\in\![m]\big\} \qquad\hbox{and}\qquad
\big\{\pi_i^*\ti\psi_0\!:~i\!\in\![m]\big\},$$
respectively. 
For example, if $m\!=\!2$,
$$\eta_3=\pi_1^*\psi_0^3+\pi_2^*\psi_0^3
+\pi_1^*\psi_0^2\,\pi_2^*\psi_0+\pi_1^*\psi_0\,\pi_2^*\psi_0^2
\in H^6\big(\ov\M_{(2,J)}(X,\be;\J)\big).$$\\

\noindent
The symmetric group on $m$ elements, $S_m$, acts on $\ov\M_{(m,J)}(X,\be;\J)$
by permuting the elements of each $m$-tuple of stable maps.
Let 
$$\cZ_{(m,J)}(X,\be;\J)=\ov\M_{(m,J)}(X,\be;\J)\big/S_m.$$
Since the map $\ev_0$ and the cohomology classes $\eta_q$ and $\ti\eta_q$
on $\ov\M_{(m,J)}(X,\be;\J)$ are $S_m$-invariant, they descendant to the quotient:
$$\ev_0\!:  \cZ_{(m,J)}(X,\be;\J)\lra X \qquad\hbox{and}\qquad
\eta_q,\ti\eta_q\in H^{2q}\big(\cZ_{(m,J)}(X,\be;\J)\big).$$\\

\noindent
The constructions of \cite{FOn} and~\cite{LT} endow
$$\ov\M_{g,k}(X,\be;\J), \qquad \ov\M_{(m,J)}(X,\be;\J), 
\qquad\hbox{and}\qquad \cZ_{(m,J)}(X,\be;\J)$$
with virtual fundamental classes (VFCs).
If the real dimension of $X$ is $2n$, the first VFC is of  real dimension
\begin{equation}\label{gkdim_e}
2\,\dim^{vir}\ov\M_{g,k}(X,\be;\J)=2\,\dim_{g,k}(X,\be)
\equiv 2\big(\blr{c_1(TX),\be}+(n\!-\!3)(1\!-\!g)+k\big).
\end{equation}
The other two VFCs are of real dimension 
\begin{equation}\label{mJdim_e}\begin{split}
2\,\dim^{vir}\ov\M_{(m,J)}(X,\be;\J)
&=2\,\dim^{vir}\cZ_{(m,J)}(X,\be;\J)\\ 
&=2\,\dim_{(m,J)}(X,\be)
\equiv 2\big(\dim_{0,|J|}(X,\be)+3-2m\big).
\end{split}\end{equation}
The VFC for $\ov\M_{(m,J)}(X,\be;\J)$ is obtained from the VFCs for
$\ov\M_{0,\{0\}\sqcup J_i}(X,\be_i;\J)$ via the usual method of 
a Kunneth decomposition of the (small) diagonal (e.g.~as in the proof of
commutativity of the quantum product).
The VFC for $\cZ_{(m,J)}(X,\be;\J)$ is the homology class induced from 
the $S_m$-action on $\ov\M_{(m,J)}(X,\be;\J)$.\\

\noindent
For each tuple
\begin{gather}\label{mutuple_e}
\mu\equiv\big(c_1,\ldots,c_k;\mu_1,\ldots,\mu_k\big)
\in (\bar\Z^+)^k\oplus H^*(X;\Z)^k\\
\label{mutuple_e2}
\hbox{s.t.}\qquad
|\mu|\equiv \sum_{j=1}^{j=k}\big(2c_i+\deg\mu_i\big)=2\,\dim_{1,k}(X,\be),
\end{gather}
let
\begin{equation*}\begin{split}
\GW_{1,k}^\be(\mu)&\equiv \blr{\big(\psi_1^{c_1}\ev_1^*\mu_1\big)
\ldots\big(\psi_k^{c_k}\ev_k^*\mu_k\big), \big[\ov\M_{1,k}(X,\be;\J)\big]^{vir}}
\qquad\hbox{and}\\
\GW_{1,k}^{\be;0}(\mu)&\equiv \blr{\big(\psi_1^{c_1}\ev_1^*\mu_1\big)
\ldots\big(\psi_k^{c_k}\ev_k^*\mu_k\big),\big[\ov\M_{1,k}^0(X,\be;\J)\big]^{vir}}
\end{split}\end{equation*}
be the standard and reduced genus-one degree-$\be$ GW-invariants of $X$ 
corresponding to the tuple~$\mu$.
The latter is constructed in~\cite{g1comp2}.
The exact definition of either invariant is not essential for the purposes
of Section~\ref{pf_sec}, as our starting point will be 
Proposition~\ref{g1comp2-bdcontr_prp} in~\cite{g1comp2},
restated as Proposition~\ref{bdcontr_prp} below,
which gives a description of the difference between the two invariants.\\

\noindent
If $\mu$ is as in~\e_ref{mutuple_e}, $m\!\in\!\Z^+$, and $J\!\subset\![k]$, let 
\begin{equation}\label{pardfn_e}
\mu_J=\prod_{j\in J}\mu_j\in H^*(X;\Z), \qquad p_J(\mu)=\sum_{j\in J}c_j, 
\qquad d_{m,J}(\mu)=n\!-\!2m-|J|+p_J(\mu).
\end{equation}
If $|\mu|\!=\!2\,\dim_{1,k}(X,\be)$, then
\begin{equation}\label{dimcount_e}
\deg\bigg(\Big(\ev_0^*\mu_J\prod_{j\not\in J}\big(\psi_j^{c_j}\ev_j^*\mu_j\big)\Big)
\cap \big[\cZ_{(m,[k]-J)}(X,\be;\J)\big]^{vir}\bigg)
=2d_{m,J}(\mu).
\end{equation}
Thus, whenever $p\!+\!q\!=\!2d_{m,J}(\mu)$, $\eta\!\in\!H^{2p}(\cZ_{(m,[k]-J)}(X,\be;\J))$,
and $\mu_0\!\in\!H^{2q}(X;\Z)$, we define
$$\GW_{(m,J)}^\be\big(\eta,\mu_0;\mu\big)=\Big\lan \eta\,\ev_0^*(\mu_0\mu_J)
\prod_{j\not\in J}\big(\psi_j^{c_j}\ev_j^*\mu_j\big), \big[\cZ_{(m,[k]-J)}(X,\be;\J)\big]^{vir}
\Big\ran\in\Q.$$

\subsection{Main Theorem}
\label{mainthm_subs}

\noindent
The description of Proposition~\ref{g1comp2-bdcontr_prp} in~\cite{g1comp2} 
implies that the difference between $\GW_{1,k}^\be(\mu)$ and 
$\GW_{1,k}^{\be;0}(\mu)$ is some linear combination
of the invariants $\GW_{(m,J)}^\be\big(\eta_p,c_q(TX);\mu\big)$ or 
equivalently of $\GW_{(m,J)}^\be\big(\ti\eta_p,c_q(TX);\mu\big)$.
The coefficients should be sums of products of top intersections of 
tautological classes on moduli spaces of genus-zero and genus-one curves, $\ov\cM_{g,N}$.
The simplest expressions in the first case, however, 
appear to be given by Hodge numbers on the blowups 
$\wt\cM_{1,([m],J)}$ of $\ov\cM_{1,m+|J|}$ constructed in 
Subsection~\ref{g1desing-curve1bl_subs} in~\cite{VaZ} and involve
the universal $\psi$-class
$$\ti\psi\equiv c_1(\wt\E)\in H^2\big(\wt\cM_{1,([m],J)}\big)$$
obtained by twisting the Hodge line bundle $\E$;
see Subsection~\ref{curvbl_subs} below.\\

\noindent
Thus, given finite sets $I$ and $J$, not both empty, and a tuple of integers 
$(\ti{c},(c_j)_{j\in J})$, we define
\begin{equation}\label{psinums_e2}
\blr{\ti{c};(c_j)_{j\in J}}_{(I,J)}=
\Big\lan \ti\psi^{\ti{c}}
\cdot\prod_{j\in J}\pi^*\psi_j^{c_j},\big[\wt\cM_{1,(I,J)}\big]\Big\ran\in\Q,
\end{equation}
where $\pi\!:\wt\cM_{1,([k],J)}\!\lra\!\ov\cM_{1,[k]\sqcup J}$ is the blow-down map.
If $\ti{c}\!+\!\sum_{j\in J}\!c_j\!\neq\!|I|\!+\!|J|$, 
$\ti{c}\!<\!0$, or $c_j\!<\!0$ for some $j\!\in\!J$,
the number in~\e_ref{psinums_e2} is zero.
By Theorem~1.1 in~\cite{g1desing2},
the numbers~\e_ref{psinums_e2} satisfy:
\begin{enumerate}[label=($R\arabic*$)]
\item\label{rec2} If $i^*\!\in\!I$ and $c_j\!>\!0$ for all $j\!\in\!J$,
$$\blr{\ti{c};(c_j)_{j\in J}}_{(I,J)}= 
\blr{\ti{c};(c_j)_{j\in J}}_{(I-i^*,J\sqcup i^*)};$$
\item\label{rec1} If $c_{j^*}\!=\!0$ for some $j^*\!\in\!J$,
$$\blr{\ti{c};(c_j)_{j\in J}}_{(I,J)}= 
|I|\blr{\ti{c}\!-\!1;(c_j)_{j\in J-\{j^*\}}}_{(I,J-j^*)}
+\!\sum_{j\in J-\{j^*\}}\!\!\!\!\!\!
      \blr{\ti{c};c_j\!-\!1,(c_{j'})_{j'\in J-\{j^*,j\}}}_{(I,J-j^*)}.$$
\end{enumerate}
Along with the relation
\begin{equation}\label{basecase_e}
\blr{0;\underset{m}{\underbrace{1,1,\ldots,1}}}_{(0,[m])}
\equiv \blr{\psi_1\ldots\psi_m,\big[\ov\cM_{1,m}\big]}=\frac{(m\!-\!1)!}{24},
\end{equation}
which follows from the usual dilaton relation (see \cite[Section 26.2]{MirSym}),
the two recursions completely determine the numbers~\e_ref{psinums_e2}.
In particular,
\begin{equation}\label{psiclass_e}
\blr{|I|\!+\!|J|;0}_{(I,J)}
\equiv\blr{\ti\psi^{|I|+|J|},\big[\wt\cM_{1,(I,J)}\big]}
=\frac{1}{24}\cdot |I|^{|J|} \cdot (|I|\!-\!1)!;
\end{equation}
see Corollary~1.2 in~\cite{g1desing2}.

\stepcounter{thm}

\begin{thmv}
\label{main_thm1}
If $(X,\om)$ is a compact symplectic manifold of real dimension $2n$,
$k\!\in\!\bar\Z^+$, $\be\!\in\!H_2(X;\Z)\!-\!0$, and 
$\mu$ is as in~\e_ref{mutuple_e} and~\e_ref{mutuple_e2}, 
then
\begin{equation}\label{main_thm1e}\begin{split}
\GW_{1,k}^{\be}(\mu)\!-\!\GW_{1,k}^{\be;0}(\mu)
&=\sum_{m\in\Z^+}\sum_{J\subset[k]}\bigg(
(-1)^{m+|J|-p_J(\mu)}\big(m\!+\!|J|\!-\!p_J(\mu);(c_j)_{j\in J}\big)_{([m],J)}\\
&\hspace{1.7in} \times \sum_{p=0}^{d_{m,J}(\mu)}\!\!\!
\GW_{(m,J)}^{\be}\big(\eta_p,c_{d_{m,J}(\mu)-p}(TX);\mu\big)\bigg).
\end{split}\end{equation}\\
\end{thmv}

\noindent
The sum in~\e_ref{main_thm1e} is finite because 
$\cZ_{(m,[k]-J)}(X,\be;\J)$ is empty if $\lr{\om,\be}/m$ is smaller than
the minimal energy of a non-constant $\J$-holomorphic map $S^2\!\lra\!X$.
Therefore, $\GW_{(m,J)}^\be$ is zero if $\lr{\om,\be}/m$ is smaller than
the minimal energy of a non-constant $\J$-holomorphic map $S^2\!\lra\!X$
for any $\om$-compatible almost complex structure~$\J$.
Theorem~\ref{main_thm1} is proved in Section~\ref{pf_sec}.\\

\noindent
We next express $\GW_{1,k}^\be(\mu)\!-\!\GW_{1,k}^{\be;0}(\mu)$ in terms of 
the numbers $\GW_{(m,J)}^\be\big(\ti\eta_p,c_q(TX);\mu\big)$.
If $J$ is a finite set, $\c\!\equiv\!(c_j)_{j\in J}$ is a $J$-tuple of integers,
and $J_0\!\subset\!J$, let
$$p_{J_0}(\c)=\sum_{j\in J_0}\!c_j\in \Z.$$
If in addition $m\!\in\!\Z^+$, we define
\begin{gather}\label{Thdfn_e}
\Th_{m,J}(\c)
=\sum_{J=\bigsqcup_{i=0}^{i=m}J_i}\!\!\!\! (-1)^{m+|J_0|-p_{J_0}(\c)}
\prod_{i=1}^{i=m}\binom{|J_i|\!-\!1}{(c_j)_{j\in J_i}}
\big(m\!+\!|J|\!-\!p_{J_0}(\c);(c_j)_{j\in J_0}\big)_{([m],J_0)},\\
\hbox{where}\qquad
\binom{|J_i|\!-\!1}{(c_j)_{j\in J_i}}\equiv
\binom{|J_i|\!-\!1}{(c_j)_{j\in J_i},|J_i|\!-\!1-p_{J_i}(\c)}
\qquad\hbox{and}\qquad
\binom{-1}{(\cdot)}\equiv1.\notag
\end{gather}
The multinomial coefficients above appear as Hodge numbers 
on $\ov\cM_{0,|J_i|+2}$.\\

\noindent
Along with the relation
$$\binom{N}{c_1,c_2,\ldots,c_l}=\binom{N\!-\!1}{c_1\!-\!1,c_2,\ldots,c_l}
+\binom{N\!-\!1}{c_1,c_2\!-\!1,c_3,\ldots,c_l}
+\!\ldots\!+\binom{N\!-\!1}{c_1,c_2,\ldots,c_{l-1},c_l\!-\!1},$$
the recursions \ref{rec2} and \ref{rec1} imply that the numbers in~\e_ref{Thdfn_e}
satisfy
\begin{enumerate}[label=($\ti{R}\arabic*$)]
\item\label{rec2b} If $m\!>\!1$ and $c_j\!>\!1$ for all $j\!\in\!J$,
$$\Th_{m,J}\big((c_j)_{j\in J}\big)=-(m\!-\!1)\Th_{m-1,J}\big((c_j)_{j\in J}\big)
+\sum_{j\in J}\Th_{m-1,J}\big(c_j\!-\!1,(c_{j'})_{j'\in J-j}\big);$$
\item\label{rec1b} If $c_{j^*}\!=\!0$ for some $j^*\!\in\!J$,
$$\Th_{m,J}\big((c_j)_{j\in J}\big)=
\!\sum_{j\in J-\{j^*\}}\!\!\!\!\!\!
\Th_{m,J-j^*}\big(c_j\!-\!1,(c_{j'})_{j'\in J-\{j^*,j\}}\big).$$
\end{enumerate}
Along with the relation
$$\Th_{1,\eset}(0)=-\frac{1}{24},$$
which follows from~\e_ref{basecase_e}, \ref{rec2b} and~\ref{rec1b}
are sufficient to determine $\Th_{m,J}\big((c_j)_{j\in J}\big)$ in many cases.
In particular,
\begin{equation}\label{somecases_e}
\Th_{m,J}\big((c_j)_{j\in J}\big)=\frac{(-1)^m(m\!-\!1)!}{24}
\begin{cases}
1,&\hbox{if}~c_j\!=\!0~\forall j;\\
0,&\hbox{if}~\sum_{j\in J}c_j<|J|.
\end{cases}\end{equation}\\

\noindent
The original recursions \ref{rec2} and \ref{rec1} are sufficient to compute 
$\Th_{m,J}\big((c_j)_{j\in J}\big)$ in all cases.
However, it is more convenient to make use of the remaining relation 
of Theorem~1.1 in~\cite{g1desing2}:
if $c_{j^*}\!=\!1$ for some $j^*\!\in\!J$, then
\begin{equation}\label{dilat_e}
\blr{\ti{c};(c_j)_{j\in J}}_{(I,J)}=\big(|I|\!+\!|J|-1\big)
\blr{\ti{c};(c_j)_{j\in J-j^*}}_{(I,J-j^*)}.
\end{equation}
This gives us a third relation for the numbers~\e_ref{Thdfn_e}:
\begin{enumerate}[label=($\ti{R}\arabic*$)]
\stepcounter{enumi}
\stepcounter{enumi}
\item\label{rec3b} If $c_{j^*}\!=\!1$ for some $j^*\!\in\!J$,
$$\Th_{m,J}\big((c_j)_{j\in J}\big)=
\big(m\!+\!|J|\!-\!1\big)\Th_{m,J-j^*}\big((c_j)_{j\in J-j^*}\big).$$
\end{enumerate}
The three relations \ref{rec2b}-\ref{rec3b}, along with the initial condition
$\Th_{1,\eset}()\!=\!-1/24$, determine the numbers 
$\Th_{m,J}(\c)$ completely.

\begin{thmv}
\label{main_thm2}
If $(X,\om)$ is a compact symplectic manifold of real dimension $2n$,
$k\!\in\!\bar\Z^+$, $\be\!\in\!H_2(X;\Z)\!-\!0$, and
$\mu$ is as in~\e_ref{mutuple_e} and~\e_ref{mutuple_e2}, 
then
\begin{equation}\label{main_thm2e}\begin{split}
\GW_{1,k}^\be(\mu)\!-\!\GW_{1,k}^{\be;0}(\mu)
&=\sum_{m\in\Z^+}\sum_{J\subset[k]}\bigg(\Th_{m,J}\big((c_j)_{j\in J}\big)\\
&\hspace{1.2in} \times
\sum_{p=0}^{d_{m,J}(\mu)}\!\!\!\GW_{(m,J)}^\be\big(\ti\eta_p,c_{d_{m,J}(\mu)-p}(TX);\mu\big)\bigg).
\end{split}\end{equation}\\
\end{thmv}

\noindent
This theorem follows immediately from Theorem~\ref{main_thm1} above
and Lemma~\ref{etadiff_lmm} below.
In turn, the latter follows from Lemma~2.2.1 and Subsection~3.2 in~\cite{P};
see also Subsection~3.3 
in~\cite{g1}.\footnote{Lemma~\ref{etadiff_lmm} is a consequence of 
the following identities. If $J_0\!\subset\!J$ is nonempty, let 
$D_{J_0}\!\subset\!\ov\M_{0,0\sqcup J}(X,\be;\J)$ be the (virtual) divisor whose
(virtually) generic element is a map from a union of two $\P^1$'s,
one of which is contracted and carries the marked points indexed 
by the set $0\!\sqcup\!J_0$.
In particular,
$$D_{J_0}\approx \ov\cM_{0,\{0,1\}\sqcup J_0}\times \ov\M_{0,0\sqcup (J-J_0)}(X,\be;\J).$$
If $\pi_P$ and $\pi_B$ are the two component projection maps, then
$$\psi_0=\ti\psi_0+\!\!\sum_{\eset\neq J_0\subset J}\!\!\!D_{J_0}, \qquad
\psi_0|_{D_{J_0}}=\pi_P^*\psi_0, 
\qquad \ti\psi_0|_{D_{J_0}}=\pi_B^*\ti\psi_0.$$}

\begin{lmm}
\label{etadiff_lmm}
Suppose $(X,\om)$ is a compact symplectic manifold of real dimension $2n$,
$k\!\in\!\bar\Z^+$, $\be\!\in\!H_2(X;\Z)\!-\!0$, and 
$\mu$ is as in~\e_ref{mutuple_e} and~\e_ref{mutuple_e2}.
If $m\!\in\!\Z^+$, $p\!\in\!\bar\Z^+$, and $J\!\subset\![k]$,
\begin{equation*}\begin{split}
\GW_{(m,J)}^\be\big(\eta_p,c_{d_{m,J}(\mu)-p}(TX);\mu\big)
&=\!\!\! \sum_{J\subset J'\subset[m]}\!\!\bigg(
\prod_{J'-J=\bigsqcup_{i=1}^{i=m}\!\!J_i}
\!\!\!\binom{|J_i|\!-\!1}{(c_j)_{j\in J_i}}\\
&\hspace{.5in} \times
\GW_{(m,J')}^\be\big(\ti\eta_{p-|J'-J|+p_{J'-J}(\mu)},c_{d_{m,J}(\mu)-p}(TX);\mu\big)
\bigg).
\end{split}\end{equation*}\\
\end{lmm}

\noindent
The advantage of Theorem~\ref{main_thm2} over Theorem~\ref{main_thm1}
is that the coefficients of the genus-zero GW-invariants in~\e_ref{main_thm2e}
satisfy simpler recursions and are more likely to vanish, due to~\e_ref{somecases_e}.
For example, if $c_j\!=\!0$ for all~$j$, i.e.~there are no descendant classes involved,
\e_ref{main_thm2e} reduces~to 
\begin{equation}\label{main_red_e}
\GW_{1,k}^\be(\mu)\!-\!\GW_{1,k}^{\be;0}(\mu)
=\frac{1}{24}\sum_{m=1}^{2m\le n}\bigg((-1)^m(m\!-\!1)!
\!\!\! \sum_{p=0}^{n-2m}\!\!
\GW_{(m,\eset)}^\be\big(\ti\eta_p,c_{n-2m-p}(TX);\mu\big)\bigg).
\end{equation}
This formula looks remarkably similar to the formula for the correction term 
in Theorem~1.1 in~\cite{g1} enumerating one-nodal rational 
curves.\footnote{In~\cite{g1}, the meanings of $\eta_p$ and $\ti\eta_p$ are reversed.}
This is not too surprising as both expressions arise from counting zeros
of analogous affine bundle maps; see Subsection~\ref{comp_subs}.

\subsection{Genus-One GW-Invariants of Calabi-Yau Hypersurfaces}
\label{cy_subs}

\noindent
The essence of mirror symmetric predictions for GW-invariants
of Calabi-Yau manifolds is that these invariants 
can be expressed in terms of certain hypergeometric series.
In this subsection, we deduce a mirror symmetry type of formula 
for the standard genus-one GW-invariants of Calabi-Yau projective 
hypersurfaces from 
a formula for the reduced genus-one GW-invariants
obtained in~\cite{bcov1}, \e_ref{main_thm2e}, and
a formula for genus-zero GW-invariants obtained in Chapter~30
of~\cite{MirSym}.
In particular, we show that the difference between the two invariants,
i.e.~\e_ref{main_thm2e}, simply cancels the last term 
in Corollary~\ref{bcov1-main_crl} in~\cite{bcov1}.
The $n\!\le\!5$ cases of Theorem~\ref{cy_thm} below
have already been obtained~\cite{bcov1},
with the $n\!=\!5$ case confirming the prediction of~\cite{BCOV}.
The $n\!=\!6$ case confirms the prediction of~\cite{KP}.\\

\noindent
Fix an integer $n\!\ge\!3$ and let $X\!\subset\!\P^{n-1}$ be a smooth
degree-$n$ hypersurface.
In this case, $d_{1,0}(X,\be)\!=\!0$ for every $\be\!\in\!H_2(X;\Z)$.
For each $d\!\in\!\Z^+$, denote by $N_1^d(X)\!\in\!\Q$ and 
$N_1^{d;0}(X)\!\in\!\Q$
the standard and reduced degree-$d$ genus-one GW-invariant of 
$X\!\subset\!\P^{n-1}$,
i.e.~the sum of $\GW_{1,0}^{\be}(\eset)$ and $\GW_{1,0}^{\be;0}(\eset)$, 
respectively, over all $\be$ lying in the preimage
of $d\ell$ under the natural homomorphism
$$H_2(X;\Z)\lra H_2(\P^{n-1};\Z),$$
where $\ell\!\in\!H_2(\P^{n-1};\Z)$ is the homology class of the line.\\

\noindent
For each $q\!=\!0,1,\ldots$, define $I_{0,q}(t)$ by 
\begin{equation}\label{Ifuncdfn_e}
\sum_{q=0}^{\i} I_{0,q}(t)w^q \equiv e^{wt}\sum_{d=0}^{\i}e^{dt}
\frac{\prod_{r=1}^{r=nd}(nw\!+\!r)}{\prod_{r=1}^{r=d}(w\!+\!r)^n}
\equiv R(w,t).
\end{equation}
Each $I_{0,q}(t)$ is a degree-$q$ polynomial in $t$ 
with coefficients that are power series in~$e^t$.
For example,
\begin{equation}\label{HGexamp_e}
I_0(t)=1+\sum_{d=1}^{\i}e^{dt}\frac{(nd)!}{(d!)^n}
\qquad\hbox{and}\qquad
I_1(t)=tI_0(t)+\sum_{d=1}^{\i}e^{dt}
\bigg(\frac{(nd)!}{(d!)^n}\sum_{r=d+1}^{nd}\!\frac{n}{r}\bigg).
\end{equation}
For $p,q\!\in\!\Z^+$ with $q\!\ge\!p$, let
\begin{equation}\label{Tden_e}
I_{p,q}(t)=\frac{d}{dt}\bigg(\frac{I_{p-1,q}(t)}{I_{p-1,p-1}(t)}\bigg).  
\end{equation}
It is straightforward to check that each of the ``diagonal'' terms $I_{p,p}(t)$ 
is a power series in $e^t$ with constant term~$1$, whenever it is defined;
see~\cite{ZaZ}, for example. 
Thus, the division in~\e_ref{Tden_e} is well-defined for all~$p$.
Proposition~\ref{bcov1-hypergeom_prp} in~\cite{bcov1} describes 
a number of relations between the power series~$I_{p,p}(t)$.
Let
\begin{equation}\label{mirmap_e}
T=\frac{I_{0,1}(t)}{I_{0,0}(t)}.
\end{equation}
By~\e_ref{HGexamp_e}, the map $t\!\lra\!T$ is a change of variables;
it will be called the {\tt mirror map}.

\begin{thm}
\label{cy_thm}
The genus-one degree-$d$ Gromov-Witten invariants of a degree-$n$ hypersurface
$X$ in $\P^{n-1}$ are given~by
\begin{equation*}\begin{split}
\sum_{d=1}^{\i}e^{dT}N_1^d(X)&=
\bigg(\frac{(n\!-\!2)(n\!+\!1)}{48}+\frac{1-(1\!-\!n)^n}{24n^2}\bigg)(T\!-\!t)
+\frac{n^2\!-\!1+(1\!-\!n)^n}{24n}\ln I_{0,0}(t)\\
&\qquad -\begin{cases}
\frac{n-1}{48}\ln\big(1\!-\!n^ne^t\big)+
\sum_{p=0}^{(n-3)/2}\frac{(n-1-2p)^2}{8}\ln I_{p,p}(t),
&\hbox{if}~2\!\not|n;\\
\frac{n-4}{48}\ln\big(1\!-\!n^ne^t\big)+
\sum_{p=0}^{(n-4)/2}\frac{(n-2p)(n-2-2p)}{8}\ln I_{p,p}(t),
&\hbox{if}~2|n,
\end{cases}
\end{split}\end{equation*}
where $t$ and $T$ are related by the mirror map~\e_ref{mirmap_e}.
\end{thm}

\noindent
The distinction between the $n$ odd and $n$ even cases appears because
the formula of Corollary~\ref{bcov1-main_crl} in~\cite{bcov1}\footnote{It 
states that $\sum_{d=1}^{\i}e^{dT}N_1^{d;0}(X)$
is given by the expression in Theorem~\ref{cy_thm} plus
$$\frac{n}{24}\sum_{p=2}^{n-2} 
\bigg(\cD_w^{n-2-p}\frac{(1\!+\!w)^n}{(1\!+\!nw)}\bigg)
\big(\cD_w^p\ln\bar{R}(w,t)\big).$$} uses the reflection symmetry
property of Proposition~\ref{bcov1-hypergeom_prp} in~\cite{bcov1}
to reduce the number of different power series~$I_{p,p}$ used.
A uniform formula can be obtained from 
Theorem~\ref{bcov1-main_thm} in~\cite{bcov1}.\\

\noindent
Let $\bar{R}(w,t)\!=\!R(w,t)/I_{0,0}(t)$. 
Then, $e^{-wt}\bar{R}(w,t)$ is a power series with $e^t$-constant term~$1$ and
$$\cD_w^p\ln\bar{R}(w,t)\equiv  \frac{1}{p!}\bigg\{\frac{d}{dw}\bigg\}^p
\Big(\ln\big(e^{-wt}\bar{R}(w,t)\big)\Big) \bigg|_{w=0} \in\Q\big[[e^t]\big]$$
for all $p\!\in\!\Z^+$ with $p\!\ge\!2$. 
Theorem~\ref{cy_thm} follows immediately from Corollary~\ref{bcov1-main_crl} in~\cite{bcov1}, 
\e_ref{main_thm2e}, and the following lemma.
Note that since $\dim\,X\!=\!n\!-\!2$ and $k\!=\!0$, 
\e_ref{main_red_e} can be written~as
\begin{equation*}
N_1^d(X)\!-\!N_1^{d;0}(X)
=\frac{1}{24}\sum_{p=2}^{n-2}\sum_{m=1}^{2m\le p}
(-1)^m(m\!-\!1)!\,
\GW_{(m,\eset)}^d\big(\eta_{p-2m},c_{n-2-p}(TX);\eset\big).
\end{equation*}

\begin{lmm}
\label{cydiff_lmm}
If $X\!\subset\!\P^{n-1}$ is a degree-$n$ hypersurface,
$x\!\in\!H^2(\P^{n-1};\Z)$ is the hyperplane class, and $p,q\!\in\!\bar\Z^+$
with $2\!\le\!p\le\!n\!-\!2$,
\begin{alignat}{1}
\label{cydiff_e1}
&c_q(TX)=\bigg(\cD_w^q\frac{(1\!+\!w)^n}{(1\!+\!nw)}\bigg)x^q|_X;\\
\label{cydiff_e2}
&\sum_{d=1}^{\i}e^{dT}\Bigg(
\sum_{m=1}^{2m\le p}(-1)^m(m\!-\!1)!\,
\GW_{(m,\eset)}^d\big(\eta_{p-2m},x^{n-2-p};\eset\big)\Bigg)
=-n\cD_w^p\ln\bar{R}(w,t)
\end{alignat}
if $T$ and $t$ are related by the mirror map~\e_ref{mirmap_e}.\\
\end{lmm}

\noindent
The relation~\e_ref{cydiff_e1} is immediate from $c(T\P^{n-1})\!=\!(1\!+\!x)^n$.
We deduce~\e_ref{cydiff_e2} below from the conclusion of 
Chapter~30 in~\cite{MirSym}.\\

\noindent
Let $\U$ be the universal curve over $\ov\M_{(m,\eset)}(\P^{n-1},d)$,
with structure map~$\pi$ and evaluation map~$\ev$:
$$\xymatrix{\U \ar[d]^{\pi} \ar[r]^{\ev} & \P^{n-1} \\
\ov\M_{(m,\eset)}(\P^{n-1},d).}$$
In other words, the fiber of $\pi$ over a tuple 
$\big([\cC_i,u_i]\big)_{i\in[m]}$ is the wedge of curves $\cC_i$
identified at the marked point~$x_0$ of each of the curves, while 
$$\ev\big(\big[(\cC_i,u_i)_{i\in[m]},z\big]\big)=u_i(z)
\qquad\hbox{if}\quad  z\!\in\!\cC_i.$$
The orbi-sheaf
$$\pi_*\ev^*\O_{\P^{n-1}}(n)\lra\ov\M_{(m,\eset)}(\P^{n-1},d)$$
is locally free, i.e.~is the sheaf of (holomorphic) sections of 
a vector orbi-bundle
$$\V_{(m,\eset)}\lra\ov\M_{(m,\eset)}(\P^{n-1},d).$$
By the (genus-zero) hyperplane-section relation,
\begin{equation}\label{HPrel_e}
(m\!-\!1)!\GW_{(m,\eset)}^d\big(\eta_{p-2m},x^{n-2-p};\eset\big)
=\frac{1}{m}\blr{\eta_{p-2m}\ev_0^*x^{n-2-p}e(\V_{(m,\eset)}),
\big[\ov\M_{(m,\eset)}(\P^{n-1},d)\big]}.
\end{equation}
The $m\!=\!1$ case of~\e_ref{HPrel_e} is completely standard,
and the same argument applies in all cases.\\

\noindent
There is a natural surjective bundle homomorphism
$$\wt\ev_0\!:\V_{(1,\eset)}\lra\ev_0^*\O_{\P^{n-1}}(n),\qquad
\big([\cC,u,\xi]\big)\lra \xi\big(x_0(\cC)\big),$$
over $\ov\M_{(1,\eset)}(\P^{n-1},d)\!\equiv\!\ov\M_{0,\{0\}}(\P^{n-1},d)$. 
Thus,
$$\V_{(1,\eset)}'\equiv\ker\wt\ev_0\lra\ov\M_{(1,\eset)}(\P^{n-1},d)$$
is a vector orbi-bundle.
It is straightforward to see~that
\begin{equation}\label{vbsplit_e}
e(\V_{(m,\eset)})=n\,\ev_0^*x\prod_{i=1}^{i=m}\pi_i^*e(\V_{(1,\eset)}').
\end{equation}
For each $r\!\in\!\bar\Z^+$, let
$$Z_r(e^T)=\sum_{d=1}^{\i}e^{dT}\blr{\psi_0^r\ev_0^*x^{n-3-r}\,e(\V_{(1,\eset)}'),
[\ov\M_{(1,\eset)}(\P^{n-1},d)]}.$$
Using the string relation (see \cite[Section~26.3]{MirSym}),
the conclusion of Chapter~30 in~\cite{MirSym} can be reformulated~as
\begin{equation}\label{g0ms_e}
e^{Tw}\bigg(1+\sum_{r=0}^{n-3}\!Z_r(e^T)w^{r+2}\bigg)
=\bar{R}(w,t)\in\Q[w]\big[\big[e^t\big]\big]/w^n,
\end{equation}
with $T$ and $t$ related by the mirror map~\e_ref{mirmap_e} as before.\\

\noindent
We now verify~\e_ref{cydiff_e2}. By~\e_ref{HPrel_e}, \e_ref{vbsplit_e}, 
and the decomposition along the small diagonal 
in $(\P^{n-1})^m$, the left-hand side of~\e_ref{cydiff_e2} equals
\begin{equation*}\begin{split}
&n\sum_{d=1}^{\i}e^{dT}\Bigg(\sum_{m=1}^{2m\le p}\frac{(-1)^m}{m}
\Blr{\ev_0^*x^{n-1-p}
\prod_{i=1}^{i=m}\pi_i^*\frac{e(\V_{(1,\eset)}')}{1\!-\!\psi_0},
\big[\ov\M_{(m,\eset)}(\P^{n-1},d)\big]}\Bigg)\\
&\qquad
=n\sum_{m=1}^{2m\le p}\frac{(-1)^m}{m}
\sum_{d=1}^{\i}e^{dT}\sum_{\underset{d_i>0}{\sum_{i=1}^{i=m}\!d_i=d}}
\sum_{\underset{p_i\ge0}{\sum_{i=1}^{i=m}\!p_i=p}}\prod_{i=1}^{i=m}
\Blr{\ev_0^*x^{n-1-p_i}\frac{e(\V_{(1,\eset)}')}{1\!-\!\psi_0},
\big[\ov\M_{(m,\eset)}(\P^{n-1},d)\big]}\\
&\qquad
=n\sum_{m=1}^{2m\le p}\frac{(-1)^m}{m}
\sum_{\underset{p_i\ge2}{\sum_{i=1}^{i=m}\!p_i=p}}
\prod_{i=1}^{i=m}\!\!Z_{p_i-2}(e^T)
=-n\,\cD^p_w\ln\bigg(1+\sum_{r=0}^{n-3}\!Z_r(e^T)w^{r+2}\bigg).
\end{split}\end{equation*}
The relation~\e_ref{cydiff_e2} now follows from~\e_ref{g0ms_e}.

\section{Proof of Theorem~\ref{main_thm1}}
\label{pf_sec}

\subsection{Blowups of Moduli Spaces of Genus-One Curves}
\label{curvbl_subs}

\noindent
In this subsection we review some aspects of the blowup construction of 
Subsection~\ref{g1desing-curve1bl_subs} in~\cite{VaZ} and add new ones,
which will be used in Subsection~\ref{comp_subs}.\\

\noindent
If $I$ is a finite set, let
\begin{equation}\label{g0and1curvsubv_e}
\A_1(I) =\big\{\big(I_P,\{I_k\!:k\!\in\!K\}\big)\!: 
K\!\neq\!\eset;~I\!=\!\bigsqcup_{k\in\{P\}\sqcup K}\!\!\!\!\!\!\!I_k;~
|I_k|\!\ge\!2 ~\forall\, k\!\in\!K\big\}.
\end{equation}
Here $P$ stands for ``principal" (component).
If $\rho\!=\!(I_P,\{I_k\!:k\!\in\!K\})$ is an element of $\A_1(I)$, 
we denote by $\cM_{1,\rho}$ the subset of $\ov\cM_{1,I}$
consisting of the stable curves~$\cC$ such~that\\
${}\quad$ (i) $\cC$ is a union of a smooth torus and $|K|$ projective lines,
indexed by~$K$;\\
${}\quad$ (ii) each line is attached directly to the torus;\\
${}\quad$ (iii) for each $k\!\in\!K$,
the marked points on the line corresponding to $k$ are indexed by $I_k$.\\
For example, the first diagram in Figure~\ref{partorder_fig}
shows an element of $\cM_{1,\rho}$ with
$$\rho=\big(\{i_1,i_2\},\big\{\{i_3,i_4,i_5,i_6\},\{i_7,i_8,i_9\}\big\}\big);$$
the number next to each component indicates the genus.
Let $\ov\cM_{1,\rho}$ be the closure of $\cM_{1,\rho}$ in~$\ov\cM_{1,I}$.
It is well-known that each space $\ov\cM_{1,\rho}$ is 
a smooth subvariety of~$\ov\cM_{1,I}$.\\

\begin{figure}
\begin{pspicture}(-1.1,-2.5)(10,1.3)
\psset{unit=.4cm}
\psarc(-2,-1){3}{-60}{60}\pscircle*(-.07,1.3){.2}\pscircle*(-.07,-3.3){.2}
\rput(.6,1.4){\smsize{$i_1$}}\rput(.6,-3.4){\smsize{$i_2$}}
\psline(-.2,0)(6.05,1.25)\pscircle*(1.8,.4){.2}\pscircle*(3.05,.65){.2}
\pscircle*(4.3,.9){.2}\pscircle*(5.55,1.15){.2}
\rput(1.8,1.1){\smsize{$i_3$}}\rput(2.9,0){\smsize{$i_4$}}
\rput(4.2,1.6){\smsize{$i_5$}}\rput(5.5,.5){\smsize{$i_6$}}
\psline(-.2,-2)(5.05,-3.25)\pscircle*(1.9,-2.5){.2}
\pscircle*(2.95,-2.75){.2}\pscircle*(4,-3){.2}
\rput(1.9,-1.8){\smsize{$i_7$}}\rput(2.95,-2.05){\smsize{$i_8$}}
\rput(4,-2.3){\smsize{$i_9$}}
\rput(-.9,-3.5){$1$}\rput(6.45,1.3){$0$}\rput(5.5,-3.2){$0$}
\rput(5,-6){\smsize{\begin{tabular}{l}$I_P\!=\!\{i_1,i_2\},K\!=\!\{1,2\}$\\
$I_1\!=\!\{i_3,i_4,i_5,i_6\},I_2\!=\!\{i_7,i_8,i_9\}$\end{tabular}}}
\rput(8,-1){\begin{Huge}$\prec$\end{Huge}}
\psarc(9,-1){3}{-60}{60}\pscircle*(10.93,1.3){.2}\pscircle*(10.93,-3.3){.2}
\rput(11.6,1.4){\begin{small}$i_1$\end{small}}
\rput(11.6,-3.4){\begin{small}$i_2$\end{small}}
\psline(10.8,0)(17.05,1.25)\pscircle*(14.05,.65){.2}
\pscircle*(15.3,.9){.2}\pscircle*(16.55,1.15){.2}
\pscircle*(12,-1){.2}\rput(12.7,-.8){\begin{small}$i_3$\end{small}}
\rput(13.9,0){\begin{small}$i_4$\end{small}}
\rput(15.2,1.6){\begin{small}$i_5$\end{small}}
\rput(16.5,.5){\begin{small}$i_6$\end{small}}
\psline(10.8,-2)(16.05,-3.25)\pscircle*(12.9,-2.5){.2}
\pscircle*(13.95,-2.75){.2}\pscircle*(15,-3){.2}
\rput(12.9,-1.8){\begin{small}$i_7$\end{small}}
\rput(13.95,-2.05){\begin{small}$i_8$\end{small}}
\rput(15,-2.3){\begin{small}$i_9$\end{small}}
\rput(10.1,-3.5){$1$}\rput(17.45,1.3){$0$}\rput(16.5,-3.2){$0$}
\psarc(19,-1){3}{-60}{60}\pscircle*(20.93,1.3){.2}\pscircle*(20.93,-3.3){.2}
\rput(21.6,1.4){\begin{small}$i_1$\end{small}}
\rput(21.6,-3.4){\begin{small}$i_2$\end{small}}
\psline(20.8,0)(27.05,1.25)\pscircle*(22.8,.4){.2}\pscircle*(24.05,.65){.2}
\pscircle*(25.3,.9){.2}\pscircle*(26.55,1.15){.2}
\rput(22.8,1.1){\begin{small}$i_3$\end{small}}
\rput(23.9,0){\begin{small}$i_4$\end{small}}
\rput(25.2,1.6){\begin{small}$i_5$\end{small}}
\rput(26.5,.5){\begin{small}$i_6$\end{small}}
\pscircle*(21.95,-.48){.2}\rput(21.4,-.7){\begin{small}$i_7$\end{small}}
\pscircle*(21.95,-1.52){.2}\rput(22.6,-1.6){\begin{small}$i_8$\end{small}}
\pscircle*(21.6,-2.5){.2}\rput(21,-2.4){\begin{small}$i_9$\end{small}}
\rput(20.1,-3.5){$1$}\rput(27.45,1.3){$0$}
\psarc(29,-1){3}{-60}{60}\pscircle*(30.93,1.3){.2}\pscircle*(30.93,-3.3){.2}
\rput(31.6,1.4){\begin{small}$i_1$\end{small}}
\rput(31.6,-3.4){\begin{small}$i_2$\end{small}}
\psline(30.8,0)(36.05,1.25)\pscircle*(32.9,.5){.2}\pscircle*(35,1){.2}
\rput(32.9,1.2){\begin{small}$i_3$\end{small}}
\rput(35,1.7){\begin{small}$i_4$\end{small}}
\psline(31,-1)(36,-1)\pscircle*(33.5,-1){.2}\pscircle*(35,-1){.2}
\rput(33.5,-.3){\begin{small}$i_5$\end{small}}
\rput(35,-.3){\begin{small}$i_6$\end{small}}
\psline(30.8,-2)(36.05,-3.25)\pscircle*(32.9,-2.5){.2}
\pscircle*(33.95,-2.75){.2}\pscircle*(35,-3){.2}
\rput(32.9,-1.8){\begin{small}$i_7$\end{small}}
\rput(33.95,-2.05){\begin{small}$i_8$\end{small}}
\rput(35,-2.3){\begin{small}$i_9$\end{small}}
\rput(30.1,-3.5){$1$}\rput(36.5,1.2){$0$}\rput(36.4,-1){$0$}\rput(36.5,-3.2){$0$}
\end{pspicture}
\caption{Examples of partial ordering~\e_ref{partorder_e}}
\label{partorder_fig}
\end{figure}

\noindent
We define a partial ordering on the set $\A_1(I)\!\sqcup\!\{(I,\eset)\}$  by setting
\begin{equation}\label{partorder_e}
\rho'\!\equiv\!\big(I_P',\{I_k'\!: k\!\in\!K'\}\big)
\prec \rho\!\equiv\!\big(I_P,\{I_k\!: k\!\in\!K\}\big)
\end{equation}
if $\rho'\!\neq\!\rho$ and there exists a map $\vph\!:K\!\lra\!K'$
such that $I_k\!\subset\!I_{\vph(k)}'$ for all $k\!\in\!K$.
This condition means that the elements of $\cM_{1,\rho'}$ can be obtained
from the elements of $\cM_{1,\rho}$ by moving more points onto the bubble components
or combining the bubble components; see Figure~\ref{partorder_fig}.\\

\noindent
Let $I$ and $J$ be finite sets such that $I$ is not empty and $|I|\!+\!|J|\!\ge\!2$. 
We put
\begin{equation*}
\A_1(I,J) =\big\{\big((I_P\!\sqcup\!J_P),\{I_k\!\sqcup\!J_k\!: k\!\in\!K\}\big)\!\in\!\A_1(I\!\sqcup\!J)\!:
~I_k\!\neq\!\eset~ \forall\, k\!\in\!K \big\}.
\end{equation*}
We note that if $\rho\!\in\!\A_1(I\!\sqcup\!J)$, then $\rho\!\in\!\A_1(I,J)$ 
if and only if every bubble component of an element of $\cM_{1,\rho}$ 
carries at least one element of~$I$.
The partially ordered set $(\A_1(I,J),\prec)$ has a unique minimal element
$$\rho_{\min}\equiv\big(\eset,\{I\!\sqcup\!J\}\big).$$
Let $<$ be an ordering on $\A_1(I,J)$ extending the partial ordering $\prec$.
We denote the corresponding maximal element by~$\rho_{\max}$.
If $\rho\!\in\!\A_1(I,J)$, define
\begin{equation}\label{minusdfn_e}
\rho\!-\!1=
\begin{cases}
\max\{\rho'\!\in\!\A_1(I,J)\!: \rho'\!<\!\rho\},&
\hbox{if}~ \rho\!\neq\!\rho_{\min};\\
0,& \hbox{if}~\rho\!=\!\rho_{\min},
\end{cases}
\end{equation}
where the maximum is taken with respect to the ordering $<$.\\

\noindent
The starting data for the blowup construction of 
Subsection~\ref{g1desing-curve1bl_subs} in~\cite{VaZ} is given~by
\begin{gather*}
\ov\cM_{1,(I,J)}^0=\ov\cM_{1,I\sqcup J},\qquad
\E_0=\E\lra \ov\cM_{1,(I,J)}^0.
\end{gather*}
Suppose $\rho\!\in\!\A_1(I,J)$ and we have constructed\\
${}\quad$ ($I1$) a blowup $\pi_{\rho-1}\!:\ov\cM_{1,(I,J)}^{\rho-1}\!\lra\!\ov\cM_{1,(I,J)}^0$
of $\ov\cM_{1,(I,J)}^0$ such that $\pi_{\rho-1}$ is one-to-one outside of the preimages
of the spaces $\ov\cM_{1,\rho'}^0$ with $\rho'\!\le\!\rho-1$;\\
${}\quad$ ($I2$) a line bundle $\E_{\rho-1}\!\lra\!\ov\cM_{1,(I,J)}^{\rho-1}$.\\
For each $\rho^*\!>\!\rho\!-\!1$, let $\ov\cM_{1,\rho^*}^{\rho-1}$ 
be the proper transform of~$\ov\cM_{1,\rho^*}$ in~$\ov\cM_{1,(I,J)}^{\rho-1}$.\\

\noindent
If $\rho\!\in\!\A_1(I,J)$ is as above, let 
$$\ti\pi_{\rho}\!:\ov\cM_{1,(I,J)}^{\rho}\lra\ov\cM_{1,(I,J)}^{\rho-1}$$ 
be the blowup of $\ov\cM_{1,(I,J)}^{\rho-1}$ along $\ov\cM_{1,\rho}^{\rho-1}$.
We denote by $\ov\cM^{\rho}_{1,\rho}$ the corresponding exceptional divisor and define
\begin{equation}\label{bundtwist_e}
\E_{\rho}=\ti\pi_{\rho}^*\,\E_{\rho-1}\otimes\O(\ov\cM^{\rho}_{1,\rho}).
\end{equation}
It is immediate that the requirements ($I1$) and ($I2$),
with $\rho\!-\!1$ replaced by~$\rho$, are satisfied.\\

\noindent
The blowup construction is concluded after $|\rho_{\max}|$ steps.
Let
$$\wt\cM_{1,(I,J)}=\ov\cM_{1,(I,J)}^{\rho_{\max}}, \qquad
\ti\E=\E_{\rho_{\max}}, \qquad \ti\psi=c_1(\ti\E).$$
By Lemma~\ref{g1desing-curve1bl_lmm} in~\cite{VaZ}, 
the end result of this blowup construction 
is well-defined, i.e.~independent of the choice of an ordering $<$ extending 
the partial ordering~$\prec$.
The reason is that different extensions of the partial order~$\prec$ correspond
to different orders of blowups along disjoint 
subvarieties.\footnote{If $\rho,\rho'\!\in\!\A_1(I,J)$ are not comparable
with respect to~$\prec$ and $\rho\!<\!\rho'$, $\ov\cM_{1,\rho}^{\rho-1}$
and $\ov\cM_{1,\rho'}^{\rho-1}$ are disjoint subvarieties in 
$\ov\cM_{1,(I,J)}^{\rho-1}$. 
However,  $\ov\cM_{1,\rho}$ and $\ov\cM_{1,\rho'}$ need not be disjoint
in $\ov\cM_{1,I\sqcup J}$.
For example, if
$$I=\{1,2,3,4\},\quad J=\eset, \quad \rho_{12}=\big((\{3,4\}),\{\{1,2\}\}\big),
\quad \rho_{34}=\big((\{1,2\}),\{\{3,4\}\}\big), \quad
\rho_{12,34}=\big((\eset),\{\{1,2\},\{3,4\}\}\big),$$
$\ov\cM_{1,\rho_{12}}$ and $\ov\cM_{1,\rho_{34}}$ intersect at 
$\ov\cM_{1,\rho_{12,34}}$ in $\ov\cM_{1,4}$, but their proper transforms
in the blowup of $\ov\cM_{1,4}$ along $\ov\cM_{1,\rho_{12,34}}$ are disjoint.}\\

\noindent
{\it Remark:} If $I\!=\!\eset$ or $|I|\!+\!|J|\!=\!1$, we define
$\wt\cM_{1,(I,J)}=\ov\cM_{1,I\sqcup J}$ and $\ti\E=\E$.\\

\noindent
We next define natural line bundle homomorphisms $s_i\!:L_i\!\lra\!\E^*$
over $\ov\cM_{1,I}$, where $L_i\!\lra\!\ov\cM_{1,I}$ is the universal tangent line
bundle at the $i$th marked point. 
These homomorphisms will then be twisted to isomorphisms~$\ti{s}_i$ on~$\wt\cM_{1,(I,J)}$.
The homomorphism~$s_i$ is induced by the natural pairing 
of tangent vectors and cotangent vectors at the $i$th marked point.
Explicitly,
\begin{gather*}
\big\{s_i([\cC;v])\big\}([\cC,\psi])=\psi_{x_i(\cC)}v \qquad\hbox{if}\\
[\cC]\!\in\!\ov\cM_{1,I}, \quad [\cC,v]\!\in\!L_i|_{\cC}\!=\!T_{x_i(\cC)}\cC, 
\quad [\cC,\psi]\!\in\!\E|_{\cC}\!=\!H^0(\cC;T^*\cC),
\end{gather*}
and $x_i(\cC)\!\in\!\cC$ is the marked point on $\cC$ labeled by~$i$.
The homomorphism~$s_i$ vanishes precisely on the curves for
which the $i$th marked point lies on a bubble component.
In fact, as divisors,
\begin{equation}\label{curvezero_e}
s_i^{-1}(0)=\!\sum_{\rho\in\B_1(I;i)}\!\!\!\!\!\ov\cM_{1,\rho},
\qquad\hbox{where}\qquad
\B_1(I;i)=\big\{ \big(I_P,\{I_B\}\big)\!\in\!\A_1(I) \!: i\!\in\!I_B\big\}.
\end{equation}\\

\noindent
If $I$ and $J$ are finite sets such that $I$ is not empty and $|I|\!+\!|J|\!\ge\!2$,
then $\B_1(I\!\sqcup\!J;i)\!\subset\!\A_1(I,J)$ for all $i\!\in\!I$.
For each $i\!\in\!I$, let
$$L_{0,i}=L_i\lra\ov\cM_{1,(I,J)}^0 \qquad\hbox{and}\qquad
s_{0,i}=s_i\in H^0\big(\ov\cM_{1,(I,J)}^0;\Hom(L_{0,i},\E_0^*)\big).$$
Suppose $\rho\!\in\!\A_1(I,J)$ and we have constructed line bundles $L_{\rho-1,i}\!\lra\!\ov\cM_{1,(I,J)}^{\rho-1}$ for $i\!\in\!I$ and 
bundle sections
\begin{equation}\label{curvezero_e2}
s_{\rho-1,i} \in
H^0(\ov\cM_{1,(I,J)}^{\rho-1};\Hom(L_{\rho-1,i},\pi_{\rho-1}^*\E^*)\big)
\quad\hbox{s.t.}\quad
s_{\rho-1,i}^{\,-1}(0)=\sum_{\rho^*\in\B_1(I\sqcup J;i),\rho^*>\rho-1}
  \!\!\!\!\!\!\!\!\!\!\!\!\!\!\! \ov\cM_{1,\rho^*}^{\rho-1}.
\end{equation}  
By~\e_ref{curvezero_e}, this assumption is satisfied for $\rho\!-\!1\!=\!0$.
If 
\begin{equation}\label{rho1dfn_e}
\rho=\big(I_P\!\sqcup\!J_P,\{I_k\!\sqcup\!J_k\!: k\!\in\!K\}\big)
\end{equation}
and $i\!\in\!I$, we define
\begin{equation}\label{bundtwist_e2}
L_{\rho,i}=\begin{cases}
\ti\pi_{\rho}^*L_{\rho-1,i}\otimes\O(\ov\cM^{\rho}_{1,\rho}),& \hbox{if}~ i\!\not\in\!I_P;\\
\ti\pi_{\rho}^*L_{\rho-1,i},& \hbox{if}~ i\!\in\!I_P.
\end{cases}
\end{equation}
By the inductive assumption, $s_{\rho-1,i}$ induces a section $s_{\rho,i}$ of
$L_{\rho,i}^*\!\otimes\!\pi_{\rho}^*\E^*$ such that
$$s_{\rho,i}^{\,-1}(0)=\sum_{\rho^*\in\B_1(I\sqcup J;i),\rho^*>\rho}
  \!\!\!\!\!\!\!\!\!\!\! \ov\cM_{1,\rho^*}^{\rho}.$$
Thus, the inductive assumption~\e_ref{curvezero_e2} is satisfied with 
$\rho\!-\!1$ replaced by~$\rho$.
Let
$$\ti{L}_i=L_{\rho_{\max},i}\lra\!\wt\cM_{1,(I,J)}, \qquad
\ti{s}_i=s_{\rho_{\max},i} \in
H^0(\wt\cM_{1,(I,J)};\Hom(\ti{L}_i,\ti\pi^*\E^*)\big).$$
By~\e_ref{curvezero_e2}, $\ti{s}_i\!:\ti{L}_i\!\lra\!\ti\pi^*\E^*$ is an isomorphism
of line bundles.\\

\noindent
{\it Remark:} The line bundles $\ti{L}_i$ and bundle isomorphisms $\ti{s}_i$ 
just defined are not the same as in Subsection~\ref{g1desing-curve1bl_subs} in~\cite{VaZ} or Subsection~2.1 in~\cite{g1desing2}.

\subsection{Blowups of Moduli Spaces of Genus-Zero Curves}
\label{curvbl0_subs}

\noindent
In this subsection we give a formula for the numbers~\e_ref{psinums_e2}
that involves the blowups of moduli spaces of genus-{\it zero} curves
defined in Subsection~\ref{g1desing-curve0bl_subs} of~\cite{VaZ}
and moduli spaces of genus-one curves, not their blowups.
The formula of Proposition~\ref{psiform_prp} will be used at the conclusion of 
Subsection~\ref{comp_subs}.\\

\noindent
If $I$ is a finite set,  let
$$\A_0(I)=\big\{\big(I_P,\{I_k\!:k\!\in\!K\}\big)\!: 
K\!\neq\!\eset;~I\!=\!\bigsqcup_{k\in\{P\}\sqcup K}\!\!\!\!\!I_k;~
|I_k|\!\ge\!2 ~\forall\, k\!\in\!K;~|K|\!+\!|I_P|\!\ge\!2\big\}.$$
Similarly to Subsection~\ref{curvbl_subs}, each element $\rho$ of $\A_0(I)$
describes a smooth subvariety
$$\ov\cM_{0,\rho}\subset \ov\cM_{0,\{0\}\sqcup I},$$
with the ``principal'' component of each curve in $\cM_{0,\rho}$ carrying
the marked points indexed by the set $\{0\}\!\sqcup\!J_P$.
There is a partial ordering~$\prec$ on $\A_0(I)$, 
defined analogously to the partial ordering~$\prec$ on~$\A_1(I)$.
If $J$ is also a finite set, let
$$\A_0(I,J)=\big\{\big((I_P\!\sqcup\!J_P),\{I_k\!\sqcup\!J_k\!: k\!\in\!K\}\big)\!\in\!\A_0(I\!\sqcup\!J)\!:
~I_k\!\neq\!\eset~ \forall\, k\!\in\!K \big\}.$$\\

\noindent
Suppose $\ale$ is a finite nonempty set and $\vr\!=\!(I_l,J_l)_{l\in\ale}$
is a tuple of finite sets such that $I_l\!\neq\!\eset$ and 
$|I_l|\!+\!|J_l|\!\ge\!2$ for all $l\!\in\!\ale$.
Let 
$$\ov\cM_{0,\vr}= \prod_{l\in\ale}\! \ov\cM_{0,\{0\}\sqcup I_l\sqcup J_l}
\qquad\hbox{and}\qquad
F_{\vr}=\bigoplus_{l\in\ale} \pi_l^*L_0 \lra \ov\cM_{0,\vr},$$
where $L_0\!\lra\!\ov\cM_{0,\{0\}\sqcup I_l\sqcup J_l}$ is 
the universal tangent line bundle for the marked point~$0$ and
$$\pi_l\!: \ov\cM_{0,\vr} \lra \ov\cM_{0,\{0\}\sqcup I_l\sqcup J_l}$$
is the projection map.
Denote by
$$\ga_{\vr}\lra\P F_{\vr}$$
the tautological line bundle.\\

\noindent
With $\vr$ as above, let
\begin{equation}\label{vr0setdfn_e}\begin{split}
\A_0(\vr) =\big\{ \big(\ale_P,(\rho_l)_{l\in\ale}\big)\!:~ 
&\ale_P\!\subset\!\ale,~\ale_P\!\neq\!\eset;~
\rho_l\!\in\!\{(I_l\!\sqcup\!J_l,\eset)\}\!\sqcup\!\A_0(I_l,J_l)~\forall\, l\!\in\!\ale;\\
&\rho_l\!=\!(I_l\!\sqcup\!J_l,\eset)~\forall\, l\!\in\!\ale\!-\!\ale_P; 
\big(\ale_P,(\rho_l)_{l\in\ale}\big)\!\neq\!
\big(\ale,(I_l\!\sqcup\!J_l,\eset)_{l\in\ale}\big)\big\}.
\end{split}\end{equation}
We define a partial ordering on $\A_0(\vr)$ by setting
\begin{equation}\label{rho0dfn_e}
\rho'\!\equiv\!\big(\ale_P',(\rho_l')_{l\in\ale}\big)
\prec  \rho\!\equiv\!\big(\ale_P,(\rho_l)_{l\in\ale}\big)
\end{equation}
if $\rho'\!\neq\!\rho$, $\ale_P'\!\subset\!\ale_P$,
and for every $l\!\in\!\ale$ either $\rho_l'\!=\!\rho_l$, $\rho_l'\!\prec\!\rho_l$,
or $\rho_l'\!=\!(I_l\!\sqcup\!J_l,\eset)$.
Let $<$ be an ordering on $\A_0(\vr)$ extending the partial ordering $\prec$.
Denote the corresponding minimal and maximal elements of $\A_0(\vr)$ 
by~$\rho_{\min}$ and~$\rho_{\max}$, respectively. 
If $\rho\!\in\!\A_0(\vr)$, define
$$\rho\!-\!1 \in \{0\} \!\sqcup\! \A_0(\vr)$$
as in~\e_ref{minusdfn_e}.\\

\noindent
If $\rho\!\in\!\A_0(\vr)$ is as in \e_ref{rho0dfn_e},
let 
$$\ov\cM_{0,\rho}=\prod_{l\in\ale}\ov\cM_{0,\rho_l}, \qquad 
F_{\rho}=\bigoplus_{l\in\ale_P} \pi_l^*L_0\big|_{\ov\cM_{0,\rho}}
\subset F_{\vr},\qquad
\wt\cM_{0,\rho}^0\!=\!\P F_{\rho} 
\subset \wt\cM_{0,\vr}^0\!\equiv\!\P F_{\vr}.$$
The spaces $\wt\cM_{0,\vr}^0$ and $\wt\cM_{0,\rho}^0$ can be represented by diagrams 
as in Figure~\ref{g0curv_fig2}.
The thinner lines represent typical elements of the spaces $\ov\cM_{0,\rho_l}$,
with the marked point~$0$ lying on the thicker vertical line.
We indicate the elements of $\ale_P\!\subset\!\ale$ with the letter~$P$
next to these points.
Note that by~\e_ref{vr0setdfn_e}, every dot on a vertical line
for which the corresponding tree has more than one line 
must be labeled with a~$P$.\\

\begin{figure}
\begin{pspicture}(-2,-2.4)(10,1.3)
\psset{unit=.4cm}
\psline[linewidth=.12](0,2.5)(0,-5.5)
\psline[linewidth=.05](0,1.5)(3.5,1.5)
\pscircle*(0,1.5){.25}\rput(-.8,1.5){\smsize{$P$}}
\pscircle*(1.5,1.5){.2}\pscircle*(3,1.5){.2}
\psline[linewidth=.05](0,-1.5)(3.5,-1.5)\pscircle*(0,-1.5){.25}\rput(-.8,-1.5){\smsize{$P$}}
\pscircle*(1,-1.5){.2}\pscircle*(2,-1.5){.2}\pscircle*(3,-1.5){.2}
\psline[linewidth=.05](0,-4.5)(3.5,-4.5)\pscircle*(0,-4.5){.25}\rput(-.8,-4.5){\smsize{$P$}}
\pscircle*(1,-4.5){.2}\pscircle*(2,-4.5){.2}\pscircle*(3,-4.5){.2}
\rput(8.5,-1.5){\begin{small}\begin{tabular}{l}
$\ale\!=\!\{1,2,3\}$\\
$|I_1\!\sqcup\!J_1|\!=\!2$\\
$|I_2\!\sqcup\!J_2|\!=\!|I_3\!\sqcup\!J_3|\!=\!3$\\ 
~\\
$\aleph_P\!=\!\{1,2,3\}$\\
$\rho_1\!=\!(I_1\!\sqcup\!J_1,\eset)$\\
$\rho_2\!=\!(I_2\!\sqcup\!J_2,\eset)$\\
$\rho_3\!=\!(I_3\!\sqcup\!J_3,\eset)$\\
\end{tabular}\end{small}}
\psline[linewidth=.1](19,2.5)(19,-5.5)
\psline[linewidth=.05](19,1.5)(22.5,1.5)\pscircle*(19,1.5){.25}
\pscircle*(20.5,1.5){.2}\pscircle*(22,1.5){.2}
\psline[linewidth=.05](19,-1.5)(22.5,-1.5)\pscircle*(19,-1.5){.25}\rput(18.2,-1.5){\smsize{$P$}}
\pscircle*(20.5,-1.5){.2}\psline[linewidth=.05](21.65,-1.85)(24,.5)
\pscircle*(22.8,-.7){.2}\pscircle*(23.6,.1){.2}
\psline[linewidth=.05](19,-4.5)(22.5,-4.5)\pscircle*(19,-4.5){.25}\rput(18.2,-4.5){\smsize{$P$}}
\pscircle*(20,-4.5){.2}\pscircle*(21,-4.5){.2}\pscircle*(22,-4.5){.2}
\rput(29,-1.5){\begin{small}\begin{tabular}{l}
$\ale\!=\!\{1,2,3\}$\\
$|I_1\!\sqcup\!J_1|\!=\!2$\\
$|I_2\!\sqcup\!J_2|\!=\!|I_3\!\sqcup\!J_3|\!=\!3$\\ 
~\\
$\aleph_P\!=\!\{2,3\}$\\
$\rho_1\!=\!(I_1\!\sqcup\!J_1,\eset)$\\
$\rho_2\!\neq\!(I_2\!\sqcup\!J_2,\eset)$\\
$\rho_3\!=\!(I_3\!\sqcup\!J_3,\eset)$\\
\end{tabular}\end{small}}
\end{pspicture}
\caption{Typical elements of $\wt\cM_{0,\vr}^0$  and $\wt\cM_{0,\rho}$}
\label{g0curv_fig2}
\end{figure}

\noindent
The blowup construction now proceeds similarly to that in 
Subsection~\ref{curvbl_subs} with
$$\E_0\!=\!\ga_{\vr} \lra \wt\cM_{0,\vr}^0.$$
The analogue of~\e_ref{bundtwist_e} has the same form:
\begin{equation}\label{bundletwist_e2b}
\E_{\rho}=\ti\pi_{\rho}^*\,\E_{\rho-1}\!\otimes\!\O(\wt\cM_{0,\rho}^{\rho}).
\end{equation}
As before, we~take
$$\wt\cM_{0,\vr}=\wt\cM_{0,\vr}^{\rho_{\max}}, \qquad
\ti\E=\E_{\rho_{\max}}, \qquad \ti\psi=c_1(\ti\E).$$
As in Subsection~\ref{curvbl_subs}, the end result of the above blowup construction
is well-defined, i.e.~independent of the choice of the ordering~$<$ 
extending the partial ordering~$\prec$.\\

\noindent
We now return to the construction of Subsection~\ref{curvbl_subs}.
If $\rho\!\in\!\A_1(I,J)$ is as in~\e_ref{rho1dfn_e}, let
$$I_P(\rho)=I_P, \quad J_P(\rho)=J_P, \quad 
\ale(\rho)=K, \quad \vr_B(\rho)=\big(I_l,J_l\big)_{l\in\ale(\rho)}.$$
Note that
\begin{equation}\label{decom_e1a}
\ov\cM_{1,\rho}\approx \ov\cM_{1,I_P(\rho)\sqcup J_P(\rho)\sqcup\ale(\rho)}
\times\ov\cM_{0,\vr_B(\rho)}.
\end{equation}
If $\rho\!=\!0$, we set
$$I_P(\rho)=I, \quad J_P(\rho)=J, \quad \ale(\rho)=\eset, \quad
\Blr{\ti\psi^{\ti{c}},\big[\wt\cM_{0,\vr_B(\rho)}]}=
\begin{cases}
1,&\hbox{if}~\ti{c}\!=\!-1;\\
0,&\hbox{otherwise}.\end{cases}$$
Let $\la\!=\!c_1(\E)$ be the Hodge class on $\ov\cM_{1,N}$.

\begin{prp}
\label{psiform_prp}
If $I$ and $J$ are finite sets and 
$(\ti{c},(c_j)_{j\in J})\!\in\!\Z\!\times\!\Z^J$, then
\begin{equation*}\begin{split}
\blr{\ti{c};(c_j)_{j\in J}}_{(I,J)}
=&\!\!\!\sum_{\rho\in\{0\}\sqcup\A_1(I,J)}\!\!\Bigg(
\Blr{\prod_{j\in J_P(\rho)}\!\!\!\!\!\psi_j^{c_j},
\big[\ov\cM_{1,I_P(\rho)\sqcup J_P(\rho)\sqcup \ale(\rho)}\big]}
\Blr{\ti\psi^{\ti{c}-1}\!\!\!\!\!\!\!\prod_{j\in J-J_P(\rho)}\!\!\!\!\!\!\!\!\psi_j^{c_j},
\big[\wt\cM_{0,\vr_B(\rho)}\big]}\\
&\qquad\qquad\qquad+
\Blr{\la\!\!\!\!\!\prod_{j\in J_P(\rho)}\!\!\!\!\!\psi_j^{c_j},
\big[\ov\cM_{1,I_P(\rho)\sqcup J_P(\rho)\sqcup \ale(\rho)}\big]}
\Blr{\ti\psi^{\ti{c}-2}\!\!\!\!\!\!\!\prod_{j\in J-J_P(\rho)}\!\!\!\!\!\!\!\!\psi_j^{c_j},
\big[\wt\cM_{0,\vr_B(\rho)}\big]}\Bigg).
\end{split}\end{equation*}\\
\end{prp}

\noindent
It is immediate that the statement holds if $\ti{c}\!\le\!0$.
For each $\rho\!\in\!\A_1(I,J)$, let
$$\wt\cM_{1,\rho}\subset\wt\cM_{1,(I,J)}$$
be the proper transform of the exceptional divisor 
$\ov\cM_{1,\rho}^{\rho}$ for the blowup at step~$\rho$.
Since $\ov\cM_{1,\rho}^{\rho}$ is a divisor in $\ov\cM_{1,(I,J)}^{\rho}$
and the blowup loci $\ov\cM_{1,\rho^*}^{\rho}$, with $\rho^*\!>\!\rho$, are
not contained in $\ov\cM_{1,\rho}^{\rho}$,
$\wt\cM_{1,\rho}$ is the pull-back of the cohomology class determined by 
$\ov\cM_{1,\rho}^{\rho}$ under the blow-down map
$$\wt\cM_{1,(I,J)}\lra \ov\cM_{1,(I,J)}^{\rho}.$$
Therefore, by~\e_ref{bundtwist_e},
\begin{equation}\label{tipsi_e}
\ti\psi=\la+\sum_{\rho\in\A_1(I,J)}\!\!\!\!\!\! \wt\cM_{1,\rho}
\in H^2\big(\wt\cM_{1,(I,J)}\big).
\end{equation}
Furthermore, by an inductive argument on the stages of the blowup construction
similar to Subsections~\ref{g1desing-map0blconstr_subs}
and~\ref{g1desing-map1blconstr_subs} in~\cite{VaZ},
\begin{equation}\label{decom_e1b}
\wt\cM_{1,\rho}\approx \wt\cM_{1,(I_P(\rho)\sqcup\ale(\rho),J_P(\rho))}
\times\wt\cM_{0,\vr_B(\rho)},\qquad
\ti\psi\big|_{\wt\cM_{1,\rho}}\approx\pi_B^*\ti\psi,
\end{equation}
where $\pi_B$ is the projection onto the second component.\footnote{The 
induction begins with
\e_ref{decom_e1a} and $\la|_{\ov\cM_{1,\rho}}\!=\!\pi_P^*\la$.
One then shows that 
$$\ov\cM_{1,\rho}^{\rho-1}\approx 
\wt\cM_{1,(I_P(\rho)\sqcup\ale(\rho),J_P(\rho))}\times\ov\cM_{0,\vr_B(\rho)}, 
\qquad 
c_1(\E_{\vr-1})|_{\ov\cM_{1,\rho}^{\rho-1}}=\pi_P^*\ti\psi,$$
and the normal bundle of $\ov\cM_{1,\rho}^{\rho-1}$ is 
$\pi_P^*\ti\E^*\!\otimes\!\pi_B^*F_{\vr_B(\rho)}$.
Thus,
$$\ov\cM_{1,\rho}^{\rho}\approx 
\wt\cM_{1,(I_P(\rho)\sqcup\ale(\rho),J_P(\rho))}
\times\wt\cM_{0,\vr_B(\rho)}^0, \qquad 
c_1(\E_{\vr})|_{\ov\cM_{1,\rho}^{\rho}}=\pi_B^*c_1(\E_0).$$
The proper transforms of $\ov\cM_{1,\rho}^{\rho}$ under blowups along
$\ov\cM_{1,\rho^*}^{\rho^*-1}$ with $\rho\!\prec\!\rho^*$ 
correspond to blowups of the second
component of~$\ov\cM_{1,\rho}^{\rho}$ and the twisting \e_ref{bundletwist_e2b}
changes $\pi_B^*c_1(\E_0)$ to~$\pi_B^*c_1(\ti\E)$.}
The $\ti{c}\!=\!1$ case of Proposition~\ref{psiform_prp} follows immediately 
from~\e_ref{tipsi_e} and~\e_ref{decom_e1b}.
If $\ti{c}\!\ge\!2$, then by~\e_ref{tipsi_e}, \e_ref{decom_e1b}, and $\la^2\!=\!0$,
\begin{equation*}\begin{split}
\ti\psi&=\ti\psi^{\ti{c}-2}\bigg(\la+\sum_{\rho\in\A_1(I,J)}\!\!\!\!\!\! \wt\cM_{1,\rho}\bigg)^2
=\ti\psi^{\ti{c}-1} \!\!\!\!\!\sum_{\rho\in\A_1(I,J)}\!\!\!\!\!\! \wt\cM_{1,\rho}
+\la\ti\psi^{\ti{c}-2}\!\!\!\!\!\sum_{\rho\in\A_1(I,J)}\!\!\!\!\!\! \wt\cM_{1,\rho}\\
&=\sum_{\rho\in\A_1(I,J)}\!\!\!\!\!\!\big(\pi_B^*\ti\psi^{\ti{c}-1}\wt\cM_{1,\rho}
+(\pi_P^*\la)(\pi_B^*\ti\psi^{\ti{c}-2})\big)\wt\cM_{1,\rho}.
\end{split}\end{equation*}
This implies the $\ti{c}\!\ge\!2$ cases of Proposition~\ref{psiform_prp}.

\subsection{Analytic Setup}
\label{diff_subs}

\noindent
We now recall the parts of Subsections~\ref{g1comp2-pertmaps_subs}, \ref{g1comp2-res_subs},
\ref{g1comp2-gwdiffsumm_subs}, and~\ref{g1comp2-setup_subs} in~\cite{g1comp2}
needed to formulate Proposition~\ref{g1comp2-bdcontr_prp} of~\cite{g1comp2}
giving a description of the difference between the two genus-one GW-invariants.\\

\noindent
An element of the moduli space $\ov\M_{g,k}(X,\be;\J)$ is represented by
a stable continuous degree-$\be$ map~$u$ from a pre-stable genus-$g$ Riemann surface $(\Si,j)$ 
with $k$ marked points to~$X$ which is smooth on each component of $\Si$ 
and satisfies the Cauchy-Riemann equation corresponding to~$(\J,j)$:
$$\bar\partial_{\J,j}u\equiv \frac{1}{2}\big(du+\J\circ du\circ j\big) =0.$$
We denote by $\ov\M_{g,k}(X,\be;\J,\nu)$ the space of solutions to the $\nu$-perturbed CR-equation:
$$\bar\partial_{\J,j}u +\nu(u)=0.$$
The perturbation term $\nu(u)$ is a section of the vector bundle
$$\La^{0,1}_{\J,j}T^*\Si\!\otimes\!u^*TX
\equiv\big\{\eta\!\in\!\hbox{Hom}_{\R}(T\Si,u^*TX)\!: \J\circ\eta=-\eta\circ j\big\}
\lra \Si$$
and depends continuously on $u$ and smoothly on each stratum $\X_{\T}(X,\be)$ 
of the space $\X_{g,k}(X,\be)$ of all continuous degree-$\be$ maps that are smooth
(or $L^p_1$ with $p\!>\!2$) on the components of the domain.
More formally, $\nu$ is a multi-section of a Banach orbi-bundle 
$\Ga^{0,1}_{g,k}(X,\be;\J)$ over $\X_{g,k}(X,\be)$.\footnote{The topological and analytic aspects
of the setup in Subsection~\ref{g1comp2-pertmaps_subs} of~\cite{g1comp2}
are analogous to~\cite{FOn} and~\cite{LT}, respectively.}\\

\noindent
Suppose $m\!\in\!\Z^+$, $J$ is a finite set, and $\be\!\in\!H_2(X;\Z)$.
If $\nu$ is a perturbation on the spaces $\X_{0,\{0\}\sqcup J_i}(X,\be_i)$
as above,
with $J_i\!\subset\!J$ and $\be_i\!\in\!H_2(X;\Z)$ such that $\om(\be_i)\!\le\!\om(\be)$, let
\begin{equation*}\begin{split}
\ov\M_{(m,J)}(X,\be;\J,\nu)&= 
\bigg\{(b_i)_{i\in[m]}\in\prod_{i=1}^{i=m}\ov\M_{0,\{0\}\sqcup J_i}(X,\be_i;\J,\nu)\!:
\be_i\!\in\!H_2(X;\Z)\!-\!\{0\},~J_i\!\subset\!J;\\
&\hspace{1.2in}~~
\sum_{i=1}^{i=m}\be_i\!=\!\be,~\bigsqcup_{i=1}^{i=m}J_i\!=\!J,~
\ev_0(b_i)\!=\!\ev_0(b_{i'})~\forall\, i,i'\!\in\![m]\bigg\}.
\end{split}\end{equation*}
Let $\M_{(m,J)}(X,\be;\J,\nu)\!\subset\!\ov\M_{(m,J)}(X,\be;\J,\nu)$ 
be the subspace consisting of $m$-tuples of maps from smooth domains.
Define
\begin{gather*}
\pi_i\!: \ov\M_{(m,J)}(X,\be;\J,\nu)\lra 
\bigsqcup_{\be_i\in H_2(X;\Z)-0}\bigsqcup_{J_i\subset J}
\ov\M_{0,\{0\}\sqcup J_i}(X,\be_i;\J,\nu),\\ 
\eta_p,\ti\eta_p\in H^{2p}\big(\ov\M_{(m,J)}(X,\be;\J,\nu)\big),\qquad
\ev_0\!: \ov\M_{(m,J)}(X,\be;\J)\lra X,
\end{gather*}
as in Subsection~\ref{gwinv_subs}.\\

\noindent
We will call a perturbation $\nu$ on $\X_{0,\{0\}\sqcup J}(X,\be)$ 
{\tt supported away from~$x_0$} if
$\nu(u)$ vanishes on a neighborhood of the marked point~$x_0$
for every element $[\Si,u]$ of $\X_{0,\{0\}\sqcup J}(X,\be)$.
In such a case, $u$ is holomorphic on a neighborhood of the marked point~$x_0$
for every element $[\Si,u]$ of $\ov\M_{0,\{0\}\sqcup J}(X,\be;\J,\nu)$.
Therefore, there is a well-defined ($\C$-linear) vector bundle homomorphism
(VBH)
$$\cD_0\!: L_0\lra\ev_0^*TX, \qquad 
\big[\Si,u;w\big]\lra du|_{x_0}w ~~~\hbox{if}~~~w\in L_0|_{[\Si,u]}=T_{x_0}\Si,$$
over $\ov\M_{0,\{0\}\sqcup J}(X,\be;\J,\nu)$.
If $m$, $J$, and $\nu$ are as in the previous paragraph and $\nu$ is
supported away from~$x_0$, we obtain $m$~VBHs
$$\pi_i^*\cD_0\!: \pi_i^*L_0\lra \ev_0^*TX$$
over $\ov\M_{0,(m,J)}(X,\be;\J,\nu)$.
The difference between the standard and reduced genus-one GW-invariants is described
below in terms of these VBHs and the homomorphisms~$s_i$
defined in Subsection~\ref{curvbl_subs}.\\

\noindent
In the genus-one case, Definition~\ref{g1comp2-pert_dfn} in~\cite{g1comp2} describes 
a class of perturbations~$\nu$ called {\tt effectively supported}.
These perturbations vanish on all components of the domain of a stable map~$u$
on which the degree of~$u$ is~zero, as well as near such components
(including after small deformations of~$[\Si,u]$).
If $\nu_{\es}$ is effectively supported, $[\Si,u]$ is an elements of
$\ov\M_{1,k}(X,\be;\J,\nu_{\es})$, and $u$ has degree~$0$ on a component~$\Si_i$ of~$\Si$,
then $u$ is constant on~$\Si$.
For a generic effectively supported perturbation~$\nu_{es}$,
$\ov\M_{1,k}(X,\be;\J,\nu_{\es})$ has the same general topological structure
as $\ov\M_{1,k}(\P^n,d)$.
In particular, if $(X,\om,\J)$ is sufficiently regular (e.g.~a low-degree projective
hypersurface), $\nu_{\es}$ can be taken to be~$0$ for our purposes.\\

\noindent
A stratum $\X_{\T}$ of $\X_{1,k}(X,\be)$ is specified by the topological type
of the domain~$\Si$ of the stable maps~$u$ in~$\X_{\T}$, including the distribution
of the $k$ marked points, and the choice of the components of~$\Si$ on which the
degree of~$u$ is not zero.
A stratum $\X_{\T}$ of $\X_{1,k}(X,\be)$ will be called {\tt degenerate} if 
the degree of any map~$u$ in $\X_{\T}$ on the principal, genus-carrying, component(s) of
its domain is zero.
The restriction of~$u$ to a component $\Si_i\!\approx\!S^2$ of $\Si$ on which 
the degree of~$u$ is not zero defines a projection
\begin{equation}\label{setup_e0}
\pi_{\T;i}\!: \X_{\T}\lra \!\!\bigsqcup_{\be_i\in H_2(X;\Z)-0}\!\!\!\!\! \!\!\!
\X_{0,K_i\sqcup J_i}(X,\be_i)
\end{equation}
for some $J_i\!\subset\![k]$ and finite set $K_i$ consisting of the nodes of~$\Si_i$.
If $\X_{\T}$ is degenerate, $K_i$ has a distinguished element, the node closest to
the principal component(s) of~$\Si$; it will be denoted by~$0$.
As in Subsection~\ref{g1comp2-setup_subs} in~\cite{g1comp2},
let $\G_{1,k}^{\gd}(X,\be;\J)$ be the space of all effectively supported deformations~$\nu_{\es}$
such that for every degenerate stratum~$\X_{\T}$, $[\Si,u]\!\in\!\X_{\T}$, and 
every component $\Si_i\!\approx\!S^2$ of $\Si$ 
on which the degree of $u$ is nonzero
\begin{equation}\label{setup_e1a}
\nu_{\es}(u)|_{\Si_i}=\big\{\pi_{\T;i}^*\nu_{\T;i}\big\}(u)|_{\Si_i}
\end{equation}
for a fixed (independent of~$u$) perturbation $\nu_{\T;i}$ on the right-hand side 
of~\e_ref{setup_e0} such that 
for every $\be_i\!\in\!H_2(X;\Z)\!-\!0$ with $\om(\be_i)\!\le\!\om(\be)$:
\begin{enumerate}[label=(gd\arabic*)]
\item\label{gdsm_item} the linearization
\begin{equation}\label{setup_e1b}
D_{\J,\nu_{\T;i};b}\!: \big\{\xi\!\in\!\Ga(\Si_b;u_b^*TX)\!: 
\xi(x_0(b))\!=\!0\big\} \lra \Ga\big(\Si_b;\La^{0,1}_{\J,j}T^*\Si_b\!\otimes\!u_b^*TX\big)
\end{equation}
of $\bpar_{\J}\!+\!\nu_{\T;i}$ at $b$ is surjective for every
$b\!\equiv\![\Si_b,u_b]\in\ov\M_{0,K_i\sqcup J_i}(X,\be_i;\J,\nu_{\T;i})$;
\item\label{gdtr_item} the restriction of the section 
$$\cD_0\in\Ga\big(\ov\M_{0,K_i\sqcup J_i}(X,\be_i;\J,\nu_{\T;i});
\Hom(L_0,\ev_0^*TX)\big)$$
to every stratum of $\ov\M_{0,K_i\sqcup J_i}(X,\be_i;\J,\nu_{\T;i})$ 
for which the degree of the maps on the component containing $x_0$ is nonzero 
is transverse to the zero section.\\
\end{enumerate}

\noindent
If $m\!\in\!\Z^+$ and $J\!\subset\![k]$, let
$$\M_{1,k}^{m,J}(X,\be;\J,\nu_{\es})\subset\ov\M_{1,k}(X,\be;\J,\nu)$$
be the subspace consisting of the stable maps $[\Si,u]$ such that $\Si$ is a union of 
a smooth torus~$\Si_P$ and $m$ spheres attached directly to~$\Si_P$, 
the degree of~$u$ is zero on~$\Si_P$ and nonzero on each of the $m$ spheres,
and $\Si_P$ carries the marked points indexed by~$J$.
If $\nu_{\es}\!\in\!\G_{1,k}^{\gd}(X,\be;\J)$, there is a natural splitting
\begin{equation}\label{g1decom_e}
\M_{1,k}^{m,J}(X,\be;\J,\nu_{\es})\approx 
\Big(\cM_{1,[m]\sqcup J}\!\times\!\M_{(m,[k]-J)}(X,\be;\J,\nu_B)\Big)\Big/S_m,
\end{equation}
where $\cM_{1,[m]\sqcup J}\!\subset\!\ov\cM_{1,[m]\sqcup J}$ is the subspace
of smooth curves and $\nu_B$ is a perturbation supported away from~$x_0$.
With our assumptions on~$\nu_{\es}$, $\nu_B$ is in fact effectively supported.
Furthermore, $\ov\M_{(m,[k]-J)}(X,\be;\J,\nu_B)$ is stratified by smooth orbifolds
(in the sense of Remark~1 in Subsection~\ref{g1comp2-pertmaps_subs} in~\cite{g1comp2})
with the expected normal bundles (i.e.~analogously to $\ov\M_{(m,[k]-J)}(\P^n,d)$).
The splitting in~\e_ref{g1decom_e} extends to an immersion over the closures,
from the right hand side to the left.\\

\noindent
Let $\mu$ be a tuple as in~\e_ref{mutuple_e} and~\e_ref{mutuple_e2}.
As in Subsection~\ref{g1comp2-gwdiffsumm_subs} in~\cite{g1comp2}, 
choose generic pseudocycle representatives $f_j\!:\bar{Y}_j\!\lra\!X$ 
for the Poincare duals of the cohomology classes~$\mu_j$ and let
$$\ov\M_{1,k}\big(X,\be;\J,\nu_{\es};(f_j)_{j\in[k]}\big)\subset 
\ov\M_{1,k}\big(X,\be;\J,\nu_{\es}\big)\times\prod_{j=1}^{j=k}\bar{Y}_j$$
be the preimage of the diagonal $\De_X^k\!\subset(X^2)^k$
under $\prod_{j=1}^{j=k}(\ev_j\!\times\!f_j)$.
Let
$$\ov\M_{1,k}\big(X,\be;\J,\nu_{\es};\mu\big)\subset 
\ov\M_{1,k}\big(X,\be;\J,\nu_{\es};(f_j)_{j\in[k]}\big)$$
be the zero set of a section $\vph$ of the vector bundle
$$V_{\mu}\equiv\bigoplus_{j=1}^{j=k}c_jL_j^*\lra 
\ov\M_{1,k}\big(X,\be;\J,\nu_{\es};(f_j)_{j\in[k]}\big).$$
For good choices of $f_j$ and $\vph$, the splitting~\e_ref{g1decom_e} induces a splitting
\begin{equation}\label{g1decom_e2}\begin{split}
\M_{1,k}^{m,J}(X,\be;\J,\nu_{\es};\mu) &\equiv 
\Big(\M_{1,k}^{m,J}(X,\be;\J,\nu_{\es})\!\times\!\prod_{j=1}^{j=k}Y_j\Big)
\cap\ov\M_{1,k}\big(X,\be;\J,\nu_{\es};\mu\big)\\
& \approx 
\Big(\cM_{1,[m]\sqcup J}(\mu)\!\times\!\M_{(m,[k]-J)}(X,\be;\J,\nu_B;\mu)\Big)\Big/S_m,
\end{split}\end{equation}
for all $m\!\in\!\Z^+$ and $J\!\subset\![k]$.
This splitting extends as an immersion over the compactifications.\\

\noindent
Denote by
$$\pi_P,\pi_B\!:
\ov\cM_{1,[m]\sqcup J}(\mu)\!\times\!\ov\M_{(m,[k]-J)}(X,\be;\J,\nu_B;\mu) \lra
\ov\cM_{1,[m]\sqcup J}(\mu)\!\times\!\ov\M_{(m,[k]-J)}(X,\be;\J,\nu_B;\mu)$$
the two component projection maps.
With $s_i$ as in Subsection~\ref{curvbl_subs} and $\cD_0$ as above, define
\begin{gather*}
\cD_{1,k}^{m,J}\!: \bigoplus_{i=1}^{i=m}
\pi_P^*L_i\!\otimes\!\pi_B^*\pi_i^*L_0\lra 
\pi_P^*\E^*\!\otimes\!\pi_B^*\ev_0^*TX,
\qquad (v_i\!\otimes\!w_i)_{i\in [m]}\lra 
\sum_{i=1}^{i=m}s_i(v)\!\otimes\!\cD_0(w_i).\notag
\end{gather*}
This is a VBH over 
$\ov\cM_{1,[m]\sqcup J}(\mu)\!\times\!\ov\M_{(m,[k]-J)}(X,\be;\J,\nu_B;\mu)$,
which descends to $\ov\M_{1,k}^{m,J}(X,\be;\J,\nu_{\es};\mu)$.\\

\noindent
Finally, suppose $\ov\cM$ is a compact topological space which is a disjoint union of 
smooth orbifolds, one of which, $\cM$, 
is a dense open subset of $\ov\cM$, and the real dimensions 
of all others do not exceed $\dim\,\cM\!-\!2$.
Let 
$$E,\O\lra \ov\cM$$
be vector orbi-bundles such that the restrictions of $E$ and $\O$ to 
every stratum of $\ov\cM$ is smooth~and
$$\rk\,\O-\rk\,E=\frac{1}{2}\dim_{\R}\cM.$$
If $\al\!\in\!\Ga\big(\ov\cM;\Hom(E,\O)\big)$
is a regular section in the sense of Definition~3.9 in~\cite{g2n2and3},
then the signed cardinality of the zero set of the affine bundle map
$$\psi_{\al,\bar\nu}\!\!\equiv\!\al\!+\!\bar\nu\!: E\lra \O$$
is finite and independent of a generic choice of $\bar\nu\!\in\!\Ga(\ov\cM;\O)$,
by Lemma~3.14 in~\cite{g2n2and3}.
We denote it by $N(\al)$.

\begin{prp}
\label{bdcontr_prp}
Suppose $(X,\om,\J)$ is a compact almost Kahler manifold of real dimension~$2n$,
$k\!\in\!\bar\Z^+$, $\be\!\in\!H_2(X;\Z)$, $\mu$ is
as in~\e_ref{mutuple_e} and~\e_ref{mutuple_e2}, and
$\nu_{\es}\!\in\!\G_{1,k}^{\gd}(X,\be;\J)$.
If the pseudocycle representatives~$f_j$   
and bundle section~$\vph$ are chosen generically, 
subject to the existence of a splitting~\e_ref{g1decom_e2}, then
\begin{equation}\label{bdcontr_prp_e1}
\GW_{1,k}^\be(\mu)-\GW_{1,k}^{\be;0}(\mu)
=\sum_{m=1}^{\i}\sum_{J\subset[k]}\cC_{1,k}^{m,J}(\bpar),
\end{equation}
where $\cC_{1,k}^{m,J}(\bpar)$ is the $\bpar$-contribution of
$\M_{1,k}^{m,J}(X,\be;\J,\nu_{\es};\mu)$ to $\GW_{1,k}^\be(\mu)$.
It is given~by
\begin{equation}\label{bdcontr_prp_e2}
\cC_{1,k}^{m,J}(\bpar)=N\big(\cD_{1,k}^{m,J}\big),
\end{equation}
where $\cD_{1,k}^{m,J}$ is a viewed as a vector bundle homomorphism
over $\ov\M_{1,k}^{m,J}(X,\be;\J,\nu_{\es};\mu)$.
In particular, $\cD_{1,k}^{m,J}$ is regular.
\end{prp}

\noindent
This is the essence of Proposition~\ref{g1comp2-bdcontr_prp} in~\cite{g1comp2}.
While Subsections~\ref{g1comp2-gwdiffsumm_subs} and~\ref{g1comp2-setup_subs}
in~\cite{g1comp2} explicitly treat 
only the case without descendants, i.e.~$c_j\!=\!0$ for all $j\!\in\![k]$,
exactly the same argument applies in the general case.
The notion of $\bpar$-contribution of a stratum to $\GW_{1,k}^\be(\mu)$
is made precise in Proposition~\ref{g1comp2-bdcontr_prp} in~\cite{g1comp2},
but~\e_ref{bdcontr_prp_e1} and~\e_ref{bdcontr_prp_e2} suffice for our purposes.

\subsection{Topological Computations}
\label{comp_subs}

\noindent
In this subsection we express the numbers~\e_ref{bdcontr_prp_e2} in terms
of cohomology classes and GW-invariants and thus conclude the proof of
Theorem~\ref{main_thm1}.\\

\noindent
With notation as in the previous subsection, let
$$\ov\M_{(m,[k]-J)}\equiv\ov\M_{(m,[k]-J)}(X,\be;\J,\nu_B;\mu).$$
Equation~\e_ref{bdcontr_prp_e2} can be restated as
\begin{equation}\label{bdcontr_prp_e3}
\cC_{1,k}^{m,J}(\bpar)=N\big(\cD_{1,k}^{m,J}\big)\big/m!
\end{equation}
with $\cD_{1,k}^{m,J}$ viewed as a VBH over 
$\ov\cM_{1,[m]\sqcup J}(\mu)\!\times\!\ov\M_{(m,[k]-J)}$.
It is straightforward to see that 
$$N(\cD_{1,k}^{m,J})=N(\wt\cD_{1,k}^{m,J}),$$
where $\wt\cD_{1,k}^{m,J}$ is the VBH over 
$\wt\cM_{1,([m],J)}(\mu)\!\times\!\ov\M_{(m,[k]-J)}$ given~by
\begin{gather*}
\wt\cD_{1,k}^{m,J}\!:  \bigoplus_{i=1}^{i=m}
\pi_P^*\ti{L}_i\!\otimes\!\pi_B^*\pi_i^*L_0
\lra \pi_P^*\E^*\!\otimes\!\pi_B^*\ev_0^*TX,
\qquad (v_i\!\otimes\!w_i)_{i\in [m]}\lra 
\sum_{i=1}^{i=m}\ti{s}_i(v)\!\otimes\!\cD_0(w_i),\notag
\end{gather*}
with $\ti{L}_i$ and $\ti{s}_i$ as at the end of 
Subsection~\ref{curvbl_subs}.\footnote{$\ov\cM_{1,[m]\sqcup J}(\mu)$ is
the zero of a section $\vph_J$ of the vector bundle 
$V_{\mu;J}\!\equiv\!\bigoplus_{j\in J}c_jL_j^*$ over 
$\ov\cM_{1,[m]\sqcup J}$ such that the restriction of $\vph_J$
to every blowup locus is transverse to the zero set.
$\wt\cM_{1,([m],J)}(\mu)$ is the preimage of $\ov\cM_{1,[m]\sqcup J}(\mu)$
under the blow-down map 
$\wt\cM_{1,([m],J)}\!\lra\!\ov\cM_{1,[m]\sqcup J}$.}\\

\noindent
Since $\ti{s}_i\!:\ti{L}_i\!\lra\!\E^*$ is an isomorphism for all $i\!\in\![m]$,
$$N(\cD_{1,k}^{m,J})=
N\big(\pi_P^*\id_{\E^*}\!\otimes\!\pi_B^*\cD_{(m,[k]-J)}\big),$$
where $\id_{\E^*}$ is viewed as a VBH on $\wt\cM_{1,([m],J)}(\mu)$
and $\cD_{(m,[k]-J)}$ is the VHB on $\ov\M_{(m,[k]-J)}$ defined~by
$$\cD_{(m,[k]-J)}\!:  F_{(m,[k]-J)}\equiv\bigoplus_{i=1}^{i=m}\pi_i^*L_0
\lra\ev_0^*TX,\qquad 
(w_i)_{i\in [m]}\lra \sum_{i=1}^{i=m}\cD_0(w_i).$$
This vector bundle homomorphism induces a VBH over the projectivization 
of~$F_{(m,[k]-J)}$:
$$\wt\cD_0\in\Ga(\wt\M_{(m,[k]-J)}^0;\Hom(\E_0;\ev_0^*TX)\big),
\quad\hbox{where}\quad
\wt\M_{(m,[k]-J)}^0=\P F_{(m,[k]-J)}, ~~~\E_0=\ga_{(m,[k]-J)},$$
and $\ga_{(m,[k]-J)}\!\lra\!\P F_{(m,[k]-J)}$ is the tautological line bundle.
It is straightforward to see from the definition that
\begin{equation}\label{projVBH_e1}
N\big(\pi_P^*\id_{\E^*}\!\otimes\!\pi_B^*\cD_{(m,[k]-J)}\big)
=N\big(\pi_P^*\id_{\E^*}\!\otimes\!\pi_B^*\wt\cD_0\big).
\end{equation}
By \ref{gdsm_item} and a dimension-count as below,
$\P F_{(m,[k]-J)}$ is stratified by orbifolds with the expected normal bundles
and $\wt\cD_0$ does not vanish on every stratum of $\P F_{(m,[k]-J)}$
on which its restriction is transverse to the zero set  
(unless $\wt\cM_{1,([m],J)}(\mu)$ is empty).
Using~\ref{gdtr_item}, the strata of~$\P F_{(m,[k]-J)}$ on which 
$\wt\cD_0$ is not transverse to the zero set can be described as follows.\\

\noindent
If $J$ is a finite set and $\be\!\in\!H_2(X;\Z)$, let
\begin{gather*}
\A_0(J)=\big\{(m;J_P,J_B)\!: m\!\in\!\Z^+;~
J_P,J_B\!\subset\!J; m\!+\!|J_P|\!\ge\!2\big\};\\
\ov\M_{(0,J)}(X,\be;\J,\nu_B)=\ov\M_{0,\{0\}\sqcup J}(X,\be;\J,\nu_B).\notag
\end{gather*}
If $\si\!=\!(m;J_P,J_B)$ is an element of $\A_0(J)$, let 
$$\M_{\si}(X,\be;\J,\nu_B) \subset \ov\M_{0,\{0\}\sqcup J}(X,\be;\J,\nu_B)$$
be the subset of  consisting of the stable maps $[\Si,u]$ such~that\\
${}\quad$ (i) the components of $\Si$ are $\Si_i\!=\!\P^1$ 
with $i\!\in\!\{P\}\!\sqcup\![k]$;\\
${}\quad$ (ii) $u|_{\Si_P}$ is constant and 
the marked points on $\Si_P$ are indexed by the set $\{0\}\!\sqcup\!J_P$;\\
${}\quad$ (iii) for each $i\!\in\![m]$, 
$\Si_i$ is attached to $\Si_P$ and $u|_{\Si_i}$ is not constant.\\
We denote by 
$$\ov\M_{\si}(X,\be;\J,\nu_B)\subset 
\ov\M_{0,\{0\}\sqcup J}(X,\be;\J,\nu_B)$$
the closure of  $\M_{\si}(X,\be;\J,\nu_B)$.
In each diagram of Figure~\ref{g0map_fig},
the irreducible components of $\Si$ are represented by lines,
and the homology class next to each component shows 
the degree of~$u$ on that component.
We indicate the marked points lying on the component $\Si_P$ only.\\

\noindent
If $m\!\in\!\Z^+$  and $J$ is a finite set, let
\begin{equation*}\begin{split}
\A_0(m;J) = \big\{ \big((\si_i)_{i\in[m]},J_B\big)\!:
\,&(\si_i,\eset)\!\in\!\{(0,\eset)\}\!\sqcup\!\A_0(J_{i,P}),~
(\si_i)_{i\in[m]}\!\neq\!(0)_{i\in[m]};
\bigsqcup_{i=1}^{i=m}J_{i,P}\!=\!J\!-\!J_B\big\}.
\end{split}\end{equation*}
If $\vr\!\equiv\!\big((\si_i)_{i\in[m]},J_B\big)$ is an element of 
$\A_0(m;J)$, we~put
$$\ale_P(\vr)=\big\{i\!\in\![m]\!: \si_i\!\neq\!0\big\}
\qquad\hbox{and}\qquad
\ale_S(\vr)=\big\{i\!\in\![m]\!: \si_i\!=\!0\big\}.$$
Here $P$ and $S$ stand for the subsets of principal and secondary 
elements of~$[m]$, respectively.
Note that $\ale_P(\vr)\!\neq\!\eset$ for all $\vr\!\in\!\A_0(m;J)$.
Let
\begin{equation*}\begin{split}
\ov\M_{\vr}(X,\be;\J,\nu_B) = 
\Big\{(b_i)_{i\in[m]}\in\prod_{i=1}^{i=m}
\ov\M_{(\si_i,J_{i,B})}(X,\be_i;\J,\nu_B)\!:
\sum_{i=1}^{i=m}\be_i\!=\!\be;~\bigsqcup_{i=1}^{i=m}J_{i,B}\!=\!J_B;\quad&\\
\ev_0(b_{i_1})\!=\!\ev_0(b_{i_2})~\forall\, i_1,i_2\!\in\![m]&\Big\}.
\end{split}\end{equation*}
This is a subspace of $\ov\M_{(m,J)}(X,\be;\J,\nu_B)$.\\

\begin{figure}
\begin{pspicture}(-1.1,-3)(10,1.5)
\psset{unit=.4cm}
\psline(3,2)(3,-4)\pscircle*(3,1.5){.2}\rput(3.6,1.6){\smsize{$j_1$}}
\pscircle*(3,-3.2){.2}\rput(3.5,-3.3){\smsize{$0$}}
\psline(1.8,0)(7.05,1.25)\rput(7.6,1.2){\smsize{$\be_1$}}
\psline(1.8,-2)(7.05,-3.25)\rput(7.6,-3.3){\smsize{$\be_2$}}
\rput(12,-1){\begin{Huge}$\prec$\end{Huge}}
\psline(18,2)(18,-4)\pscircle*(18,1.5){.2}\rput(18.6,1.6){\smsize{$j_1$}}
\pscircle*(18,-3.2){.2}\rput(18.5,-3.3){\smsize{$0$}}
\psline(16.8,0)(22.05,1.25)\rput(22.6,1.2){\smsize{$\be_1'$}}
\psline(16.8,-2)(22.05,-3.25)\rput(22.6,-3.3){\smsize{$\be_2'$}}
\rput(22,-6){\smsize{\begin{tabular}{l}$m\!=\!2,\,J_P\!=\!\{j_1\}$\\ 
$\be_1',\be_2'\!\neq\!0$, $\be_1'\!+\!\be_2'\!=\!\be$\end{tabular}}}
\psline(28,2)(28,-4)
\pscircle*(28,-3.2){.2}\rput(28.5,-3.3){\smsize{$0$}}
\psline(26.8,0)(32.05,1.25)\rput(32.6,1.2){\smsize{$\be_1''$}}
\psline(27,-1)(32,-1)\rput(32.5,-1.1){\smsize{$\be_2''$}}
\psline(26.8,-2)(32.05,-3.25)\rput(32.6,-3.3){\smsize{$\be_3''$}}
\end{pspicture}
\caption{Examples of partial ordering~\e_ref{partorder_e2}}
\label{g0map_fig}
\end{figure}

\noindent
With $\mu$ as before and $\vr\!\in\!\A_0(m;[k]\!-\!J)$, let
$$\ov\M_{\vr}
=\Big(\ov\M_{\vr}(X,\be;\J,\nu_B)\!\times\!\prod_{j=1}^{j=k}\bar{Y}_j\Big)
\cap \ov\M_{(m,[k]-J)}(X,\be;\J,\nu_B;\mu).$$
Define
$$F_{\vr;P} =\bigoplus_{i\in\ale_P(\vr)}
\!\!\! \pi_i^*L_0\Big|_{\ov\M_{\vr}} 
\subset F_{(m,[k]-J)}\big|_{\ov\M_{\vr}},
\qquad
\wt\M_{\vr}^0=\P F_{\vr;P}\subset \wt\M_{(m,[k]-J)}^0.$$
It is immediate from the definition of~$\wt\cD_0$ that it vanishes
identically on $\wt\M_{\vr}^0$ for every element $\vr$ in $\A_0(m;[k]\!-\!J)$,
since $\cD_0$ vanishes identically on the strata of 
$\ov\M_{0,\{0\}\sqcup J_i}(X,\be_i;\J,\nu_B)$ for which the degree
of the maps on the component carrying the $0$th marked point is zero.
On the other hand, by~\ref{gdtr_item}, the restriction of $\wt\cD_0$ 
to any stratum of  $\wt\M_{(m,[k]-J)}^0$ in the complement of 
every $\wt\M_{\vr}^0$ is transverse to the zero set and thus does not vanish
by a dimension count as below (unless $\wt\cM_{1,([m],J)}(\mu)$ is empty).\\

\noindent
As described in Section~3 of~\cite{g2n2and3}, the number $N(\wt\cD_0)$
is the euler class of the quotient of the target bundle of $\wt\cD_0$
by the domain line bundle minus a correction from~$\wt\cD_0^{-1}(0)$.
The correction splits into contributions from the strata of~$\wt\cD_0^{-1}(0)$
each of which is again the number of zeros of an affine bundle map,
but with the rank of the target bundle reduced by at least one.
The linear part of each affine bundle map is determined by
the topological behavior of~$\wt\cD_0^{-1}$ in the normal direction
to each stratum.
This behavior (for~$\cD_0$ and thus for~$\wt\cD_0$) is described 
in Subsection~\ref{g1comp2-g0prp_subs} of~\cite{g1comp2}. 
Thus, by iteration, one obtains a finite tree of cohomology classes at
the nodes which sum up to~$N(\wt\cD_0)$.
The tree in this case is similar to a subtree of the tree 
in Subsection~3.2 of~\cite{g1}, but twisted with~$\E^*$.
Thus, $N(\wt\cD_0)$ can be expressed in terms of cohomology classes
by a direct, though laborious, computation nearly identical to
the one in Subsections~3.1 and~3.2 in~\cite{g1}.
This time, we will instead compute $N(\wt\cD_0)$ 
by blowing up $\wt\M_{(m,[k]-J)}^0$ and twisting $\wt\cD_0$ to 
a nowhere-vanishing vector bundle homomorphism $\wt\cD_{(m,[k]-J)}$.
This construction is a direct generalization of 
Section~\ref{g1desing-map0bl_sec} in~\cite{VaZ}.\\

\noindent 
Define a partial ordering on the set $\A_0(J)$ by setting
\begin{equation}\label{partorder_e2}
\si'\!\equiv\!(m';J_P',J_B') \prec \si\!\equiv\!(m;J_P,J_B)
\qquad\hbox{if}\quad \si'\!\neq\!\si,~m'\!\le\!m,~J_P'\!\subset\!J_P.
\end{equation}
Similarly to Subsection~\ref{curvbl_subs}, 
this condition means that the elements of $\M_{\si'}(X,\be;\J,\nu)$ 
can be obtained from the elements of $\M_{\si}(X,\be;\J,\nu_B)$ 
by moving more points 
onto the bubble components or combining the bubble components; 
see Figure~\ref{g0map_fig}.
The bubble components are the components not containing the marked point~$0$.
Define a partial ordering $\prec$ on $\A_0(m;J)$ by setting
\begin{equation}\label{vrdfn_e2}
\vr'\!\equiv\!\big((\si_i')_{i\in[m]},J_B'\big) 
\prec \vr\!\equiv\!\big((\si_i)_{i\in[m]},J_B)
\end{equation}
if $\vr'\!\neq\!\vr$ and for every $i\!\in\![m]$ either 
$\si_i'\!=\!\si_i$, $(\si_i',\eset)\!\prec\!(\si_i,\eset)$, or
$\si_i'\!=\!0$. Note that 
\begin{equation}\label{vrdfn_e3}
\vr'\prec \vr ~~~\Lra~~~ \ale_P(\vr')\subset\ale_P(\vr); \qquad
\vr=\big((m_i;J_{i,P})_{i\in\ale_P(\vr)},(0)_{i\in\ale_S(\vr)},J_B\big)
\end{equation}
for some $m_i$ and $J_{i,P}$.
Choose an ordering $<$ on $\A_0(m;J)$ extending the partial ordering~$\prec$.
Denote the corresponding minimal element by~$\vr_{\min}$
and the largest element for which $\ov\M_{\vr}$ is nonempty by~$\vr_{\max}$.
For every  $\vr\!\in\!\A_0(m;J)$, define
$$\vr\!-\!1 \in \{0\}\!\sqcup\!\A_0(m;J)$$
as in~\e_ref{minusdfn_e}.\\

\noindent
With $\vr$ as \e_ref{vrdfn_e3}, let
$$\vr_P=\big([m_i],J_{i,P}\big)_{i\in\ale_P(\vr)}, \quad
m_B(\vr)=\big|\ale_S(\vr)|+\sum_{i\in\ale_P(\vr)}\!\!\!\!m_i, \quad
J_B(\vr)=J_B, \quad\hbox{and}\quad 
G_{\vr}=\prod_{i\in\ale_P(\vr)}\!\!\!\!\!S_{m_i}.$$
With $\wt\cM_{0,\vr_P}^0$ as in Subsection~\ref{curvbl0_subs},
there is a natural node-identifying immersion
$$\io_{0,\vr}\!: \wt\cM_{0,\vr_P}^0(\mu) \times 
\ov\M_{(m_B(\vr),J_B(\vr))} 
\lra \wt\M_{\vr}^0 \subset \wt\M_{(m,[k]-J)}^0.$$
It descends to an immersion
$$\bar\io_{0,\vr}\!: 
\big(\wt\cM_{0,\vr_P}^0(\mu) \!\times\! 
\ov\M_{(m_B(\vr),J_B(\vr))}\big)\big/G_{\vr}
\lra \wt\M_{(m,[k]-J)}^0,$$
which is an embedding outside the preimages of $\wt\M_{\vr'}^0$ 
with $\vr'\!\prec\!\vr$.\\

\noindent
As in the blowup construction of Subsection~\ref{curvbl_subs},
we inductively define 
$$\ti\pi_{\vr}\!: \wt\M_{(m,[k]-J)}^{\vr}
\lra\wt\M_{(m,[k]-J)}^{\vr-1}$$
to be the blowup of $\wt\M_{(m,[k]-J)}^{\vr-1}$ along the proper
transform $\wt\M_{\vr}^{\vr-1}$ of $\wt\M_{\vr}^0$
in $\wt\M_{(m,[k]-J)}^{\vr-1}$.
If $\wt\M_{\vr}^{\vr}\!\subset\!\wt\M_{(m,[k]-J)}^{\vr}$
is the exceptional divisor, let
\begin{equation}\label{bndltwist_e4}
\E_{\vr}=\ti\pi_{\vr}^*\E_{\vr-1}\otimes \O\big(\wt\M_{\vr}^{\vr}\big).
\end{equation}
The vector bundle homomorphism $\wt\cD_{\vr-1}\!:\E_{\vr-1}\!\lra\!\ev_0^*TX$
induces a section
$$\wt\cD_{\vr}\in\Ga\big(\wt\M_{(m,[k]-J)}^{\vr};
\Hom(\E_{\vr},\ev_0^*TX)\big).$$\\

\noindent
As described in detail in Subsection~\ref{g1desing-map0blconstr_subs} in~\cite{VaZ} 
(in the case $(X,\J)\!=\!\P^n$), $\io_{0,\vr}$ induces an immersion
$$\io_{\vr-1,\vr}\!: \wt\cM_{0,\vr_P}(\mu) \times  \ov\M_{(m_B(\vr),J_B(\vr))} 
\lra \wt\M_{\vr}^{\vr-1} \subset \wt\M_{(m,[k]-J)}^{\vr-1}$$
and an {\it embedding}
$$\bar\io_{\vr-1,\vr}\!: 
\big(\wt\cM_{0,\vr_P}(\mu)\!\times\!\ov\M_{(m_B(\vr),J_B(\vr))}\big)/G_{\vr}
\lra \wt\M_{(m,[k]-J)}^{\vr-1}.$$
Thus, the centers of all blowups are smooth (in the appropriate sense)
and 
$$\wt\M_{\vr}^{\vr}\approx 
\big(\wt\cM_{0,\vr_P}(\mu) \times\wt\M_{(m_B(\vr),J_B(\vr))}^0\big)\big/G_{\vr}.$$
Furthermore, 
\begin{equation}\label{psirestr_e4}
\ti\pi_{\vr}^*c_1(\E_{\vr-1})\big|_{\wt\M_{\vr}^{\vr}}=\pi_P^*\ti\psi,
\qquad 
c_1(\E_{\vr})\big|_{\wt\M_{\vr}^{\vr}}
=\pi_B^*c_1\big(\ga_{(m_B(\vr),J_B(\vr))}\big),
\end{equation}
where
$$\pi_P,\pi_B\!: \wt\cM_{0,\vr_P}(\mu)\times  \wt\M_{(m_B(\vr),J_B(\vr))}^0 
\lra \wt\cM_{0,\vr_P}(\mu),\wt\M_{(m_B(\vr),J_B(\vr))}^0$$
are the projection maps.
Finally, the restriction of $\cD_{\vr}$ to every stratum of 
$\wt\M_{(m,[k]-J)}^{\vr}$ not contained in the proper transform
$\wt\M_{\vr^*}^{\vr}$ of $\wt\M_{\vr}^0$ for any 
$\vr^*\!\in\!\A_0(m;[k]\!-\!J)$ with $\vr^*\!>\!\vr$
is transverse to the zero~set.\footnote{This statement is obtained as
in Subsection~\ref{g1desing-map0blconstr_subs} in~\cite{VaZ},
using the description of the behavior of $\cD_0$ in 
Subsection~\ref{g1comp2-g0prp_subs} in~\cite{g1comp2} and 
the assumption~\ref{gdtr_item}.}\\

\noindent
Define
$$\wt\M_{(m,[k]-J)}=\wt\M_{(m,[k]-J)}^{\vr_{\max}}, \quad
\wt\E=\E_{\vr_{\max}}, \quad
\wt\cD_{(m,[k]-J)}\in
\Ga\big(\wt\M_{(m,[k]-J)};\Hom(\wt\E_{(m,[k]-J)},\ev_0^*TX)\big).$$
As can be seen directly from the definition,
$$N(\pi_P^*\id_{\E^*}\!\otimes\!\pi_B^*\wt\cD_0\big)
=N\big(\pi_P^*\id_{\E^*}\!\otimes\!\pi_B^*\wt\cD_{(m,[k]-J)}\big),$$
where the maps $\pi_P$ and $\pi_B$ on the right-hand side are 
the two component projections
$$\wt\cM_{1,([m],J)}(\mu)\times\wt\M_{(m,[k]-J)}\lra
\wt\cM_{1,([m],J)}(\mu),\wt\M_{(m,[k]-J)}.$$
On the other hand, by the previous paragraph, 
the restriction of $\wt\cD_{(m,[k]-J)}$ to every stratum of 
$\wt\M_{(m,[k]-J)}$ is transverse to the zero set.
By~\e_ref{gkdim_e}, \e_ref{mJdim_e}, \e_ref{mutuple_e2}, and~\ref{gdsm_item},
$$\dim\, \wt\cM_{1,([m],J)}(\mu)\!\times\!\wt\M_{(m,[k]-J)}
=2(n\!-\!m)+2(m\!-\!1)
<2\rk_{\C}\big(\pi_P^*\E^*\!\otimes\!\pi_B^*\ev_0^*TX\big).$$
Therefore, $\pi_P^*\id_{\E^*}\!\otimes\!\pi_B^*\wt\cD_{(m,[k]-J)}$ 
is injective and thus
\begin{equation}\label{cohomform_e1}\begin{split}
\cC_{1,k}^{m,J}(\bpar)&=\frac{1}{m!}
N\big(\pi_P^*\id_{\E^*}\!\otimes\!\pi_B^*\wt\cD_{(m,[k]-J)}\big)\\
&=\frac{1}{m!}\blr{e\big((\pi_P^*\E^*\!\otimes\!\pi_B^*\ev_0^*TX)\big/
(\pi_P^*\E^*\!\otimes\!\pi_B^*\wt\E)\big),
\big[\wt\cM_{1,([m],J)}(\mu)\!\times\!\wt\M_{(m,[k]-J)}\big]}.
\end{split}\end{equation}\\

\noindent
{\it Remark:} Since $c_1(\E^*)^2\!=\!\la^2\!=\!0$, the last expression 
in~\e_ref{cohomform_e1} is zero if $m\!+\!|J|\!>\!1$ and $c_j\!=\!0$
for all $j\!\in\!J$.
Thus, if $\mu$ involves no descendants, i.e.~$c_j\!=\!0$ for 
all $j\!\in\![k]$, the only stratum of $\ov\M_{1,k}(X,\be;\J,\nu_{\es})$
contributing to the difference between 
the standard and reduced genus-one degree-$\be$ invariants
corresponding to~$\mu$ is $\M_{1,k}^{1,\eset}(X,\be;\J,\nu_{\es})$.\\

\noindent
It remains to express the right-hand side of~\e_ref{cohomform_e1} 
in terms of GW-invariants.
Let  
\begin{gather*}
\blr{\mu}_{1,[m]\sqcup J}=
\Blr{\prod_{j\in J}\!\psi_j^{c_j},\big[\ov\cM_{1,[m]\sqcup J}]}, \qquad
\blr{\la;\mu}_{1,[m]\sqcup J}=
\Blr{\la\prod_{j\in J}\!\psi_j^{c_j},\big[\ov\cM_{1,[m]\sqcup J}]},\\
\blr{\ti\psi^p;\mu}_{(0,\vr_P)}=
\blr{\ti\psi^p\!\!\!\!\prod_{j\in J_P(\vr)}\!\!\!\!\!\psi_j^{c_j},
\big[\wt\cM_{0,\vr_P}\big]}, \\
\bar\A_0(m;[k]\!-\!J)=\{0\}\sqcup\A_0(m;[k]\!-\!J),
\qquad \bar\A_1(I,J)=\{0\}\sqcup\A_1(I,J).
\end{gather*}
Since $\la^2\!=\!0$, by~\e_ref{cohomform_e1},
\begin{equation}\label{cohomform_e3}\begin{split}
\cC_{1,k}^{m,J}(\bpar)
&=\frac{1}{m!}\bigg(\blr{\mu}_{1,[m]\sqcup J}
\sum_{p=1}^n \blr{c_1(\wt\E^*)^{p-1} \ev_0^*c_{n-p}(TX),
\big[\wt\M_{(m,[k]-J)}]}\\
&\qquad\qquad- \blr{\la;\mu}_{1,[m]\sqcup J}
\sum_{p=1}^{n-1} \blr{c_1(\wt\E^*)^{p-1} \ev_0^*c_{n-1-p}(TX),
\big[\wt\M_{(m,[k]-J)}]}\bigg).
\end{split}\end{equation}
For each $\vr\!\in\!\bar\A_0(m;[k]\!-\!J)$,
let $\bar\pi_{\vr}\!:\wt\M_{(m,[k]-J)}\!\lra\!\wt\M_{(m,[k]-J)}^{\vr}$
be the blow-down map.
By~\e_ref{bndltwist_e4},
\begin{gather*}
c_1(\wt\E)=\bar\pi_0^*c_1\big(\ga_{(m,[k]-J)}\big)+\!
\sum_{\vr\in\A_0(m;[k]\!-\!J)}\!\!\!\!\!\!\!\!\!\!
\bar\pi_{\vr}^*\wt\M_{\vr}^{\vr}\qquad\Lra\\
\begin{split}
c_1(\wt\E)^{p-1}&=\bar\pi_0^*c_1\big(\ga_{(m,[k]-J)}\big)^{p-1}
+\!\!\!\sum_{\vr\in\A_0(m;[k]\!-\!J)}\!
\sum_{q=1}^{p-1}
\bigg(\la\!+\!\!\sum_{\vr'<\vr}\!\bar\pi_{\vr'}^*\wt\M_{\vr'}^{\vr'}\bigg)^{p-1-q}
\!\bigg(\la\!+\!\!\sum_{\vr'\le\vr}\!\bar\pi_{\vr'}^*\wt\M_{\vr}^{\vr}\bigg)^{q-1}
\!\!\!\bar\pi_{\vr}^*\wt\M_{\vr}^{\vr}\\
&=\bar\pi_0^*c_1\big(\ga_{(m,[k]-J)}\big)^{p-1}
+\sum_{\vr\in\A_0(m;[k]\!-\!J)}\!\bar\pi_{\vr}^*
\Bigg(\sum_{q=1}^{p-1}
\big(\ti\pi_{\vr}^*c_1(\E_{\vr-1})^{p-1-q}c_1(\E_{\vr})^{q-1}\big) 
\cap\wt\M_{\vr}^{\vr}\Bigg).
\end{split}\end{gather*}
Note that for every $\vr\!\in\!\bar\A_0(m;[k]\!-\!J)$
\begin{equation*}\begin{split}
&\blr{c_1(\ga_{(m_B(\vr),J_B(\vr))}^*)^{q-1}\ev_0^*c_r(TX),
\big[\wt\M_{(m_B(\vr),J_B(\vr))}^0\big]}
=\blr{\eta_{q-m_B(\vr)}\ev_0^*c_r(TX),
\big[\ov\M_{(m_B(\vr),J_B(\vr))}\big]}\\
&\hspace{2.5in}
=m!|G_{\vr}|\,\GW_{(m_B(\vr),J_B(\vr))}^{\be}
\big(\eta_{q-m_B(\vr)},c_r(TX);\mu\big),
\end{split}\end{equation*}
with $(m_B(0),J_B(0))\!\equiv\!(m;[k]\!-\!J)$ and $|G_0|\!\equiv\!1$.\\

\noindent
Thus, by~\e_ref{psirestr_e4} and~\e_ref{cohomform_e3},
\begin{alignat}{1}\label{cohomform_e4}
\cC_{1,k}^{m,J}(\bpar)&=
\sum_{\rho\in\bar\A_0(m;[k]-J)}\sum_{p=1}^{p=n}\sum_{q=1}^{q=p}
\Bigg\{(-1)^{p-q}\blr{\ti\psi^{p-1-q};\mu}_{(0,\vr_P)}\\
&\hspace{1.2in}
\times\Bigg(\blr{\mu}_{1,[m]\sqcup J}
\GW_{(m_B(\vr),J_B(\vr))}^{\be}
\big(\eta_{q-m_B(\vr)},c_{n-p}(TX);\mu\big)\notag\\
&\hspace{1.5in}
-\blr{\la;\mu}_{1,[m]\sqcup J}
\GW_{(m_B(\vr),J_B(\vr))}^{\be}
\big(\eta_{q-m_B(\vr)},c_{n-1-p}(TX);\mu\big)\Bigg)\Bigg\},\notag
\end{alignat}
where we set
$$\blr{\ti\psi^r;\mu}_{(0,0_P)}
=\begin{cases} 
1,&\hbox{if}~r\!=\!-1;\\
0,&\hbox{otherwise}.
\end{cases}$$
Most terms in~\e_ref{cohomform_e4} vanish for dimensional reasons.
By~\e_ref{dimcount_e},
\begin{alignat*}{1}
&\cC_{1,k}^{m,J}(\bpar)=\sum_{m^*=m}^{\i}\sum_{J\subset J^*\subset[k]}
\sum_{\underset{\stackrel{|I_P(\rho)|+|\ale(\rho)|=m}
{J_P(\rho)=J}}{\rho\in\bar\A_1([m^*],J^*)}}
\sum_{q=0}^{n-m^*}\bigg\{(-1)^{n-m^*-d_{m^*,J^*}(\mu)}
\GW_{(m^*,J^*)}^{\be}\big(\eta_q,c_{d_{m^*,J^*}(\mu)-q};\mu\big)\\
&~~\times
\bigg(\blr{\mu}_{1,[m]\sqcup J}
\blr{\ti\psi^{n-m^*-d_{m^*,J^*}(\mu)-1};\mu}_{(0,\vr_B(\rho))}
+\blr{\la;\mu}_{1,[m]\sqcup J}
\blr{\ti\psi^{n-m^*-d_{m^*,J^*}(\mu)-2};\mu}_{(0,\vr_B(\rho))}
\bigg)\Bigg\}.
\end{alignat*}
Summing over all $(m,J)$ as required by~\e_ref{bdcontr_prp_e1} and 
using the last expression in~\e_ref{pardfn_e}, we obtain
\begin{equation*}\begin{split}
\GW_{1,k}^\be(\mu)-\GW_{1,k}^{\be;0}(\mu)
&=\sum_{m^*=1}^{\i}\sum_{J^*\subset[k]}\Bigg\{
(-1)^{m^*+|J^*|-p_{J^*}(\mu)}\\
&\times\bigg( \sum_{q=0}^{d_{m^*,J^*}(\mu)}\!\!\!\!
\GW_{(m^*,J^*)}^{\be}\big(\eta_q,c_{d_{m^*,J^*}(\mu)-q};\mu\big)\bigg)\\
&\times\sum_{\rho\in\bar\A_1([m^*],J^*)}\!\!
\bigg(\blr{\mu}_{1,I_P(\rho)\sqcup J_P(\rho)\sqcup \aleph(\rho)}
\blr{\ti\psi^{m^*+|J^*|-p_{J^*}(\mu)-1};\mu}_{(0,\vr_B(\rho))}\\
&\hspace{.8in}
+\blr{\la;\mu}_{1,I_P(\rho)\sqcup J_P(\rho)\sqcup \aleph(\rho)}
\blr{\ti\psi^{m^*+|J^*|-p_{J^*}(\mu)-2};\mu}_{(0,\vr_B(\rho))}
\bigg)\Bigg\}.
\end{split}\end{equation*}
Finally, Proposition~\ref{psiform_prp} reduces the last expression to
the statement of Theorem~\ref{main_thm1}.\\

\noindent
{\it Remark:} Since $\wt\M_{(m,[k]-J)}^0$ is not a complex manifold,
some care is needed in constructing its ``complex'' blowups.
These are obtained by modifying normal neighborhoods to 
the strata of the blowup centers in the expected way.
The information needed to specify the normal bundles to such strata
is described in Subsection~\ref{g1comp2-g0prp_subs} of~\cite{g1comp2}.
Similarly, \e_ref{bndltwist_e4} describes a twisting of line bundles,
not of sheaves.
In fact, we know a priori that $N(\wt\cD_0)$ depends only on
the topology of the situation:
\begin{enumerate}[label=(T\arabic*)]
\item\label{al_item} the domain and target bundles of $\wt\cD_0$;
\item\label{nb_item} the normal bundles to the strata of~$\wt\cD_0^{-1}(0)$;
\item\label{be_item} the topological behavior of $\wt\cD_0$ in the normal
directions to the strata of~$\wt\cD_0^{-1}(0)$.
\end{enumerate}
By constructing a tree of chern classes, as suggested above and
similarly to Subsection~3.2 in~\cite{g1}, one can obtain 
a universal formula expressing $N(\wt\cD_0)$ in terms of 
the chern class of~\ref{al_item} and~\ref{nb_item} evaluated on the closures of
the strata of~$\wt\cD_0^{-1}(0)$, with the coefficients determined by~\ref{be_item}.
If such a universal formula holds in the presence of additional geometry
(e.g.~in the complex category), it must hold in general.
Thus, it is sufficient to obtain a formula for $N(\wt\cD_0)$
assuming $\wt\M_{(m,[k]-J)}^0$ is a complex manifold.\\

\vspace{.2in}

\noindent
{\it Department of Mathematics, SUNY, Stony Brook, NY 11794-3651}\\
azinger@math.sunysb.edu\\

\end{document}